\documentclass[11pt,a4paper]{article}
\usepackage[top=2cm,left=2.3cm, right=2.3cm, bottom=2cm]{geometry}
\usepackage[english]{babel}
\usepackage[utf8]{inputenc}
\usepackage[T1]{fontenc}
\usepackage{microtype}
\usepackage{lmodern}
\usepackage{marvosym}
\usepackage{amsmath,amssymb,amsthm,amstext,amsfonts}
\usepackage{mathtools}
\usepackage{mathrsfs}
\usepackage[shortlabels]{enumitem}
\usepackage{xspace}
\usepackage{tikz}
\usepackage{helvet}
\usepackage{graphicx}
\usepackage{float}
\usepackage{subfigure}
\usepackage{listings}
\usepackage{dsfont}
\usepackage{setspace}
\usepackage{color}
\usepackage{hyperref}
\usepackage{emptypage}
\usepackage{tabularx}
\usepackage{pbox}
\usepackage{ marvosym } 
\usepackage[nooneline, bf]{caption}	
\usepackage{comment}



\newcommand*{\cA}{\mathcal{A}}

\newcommand*{\R}{\mathbb{R}}

\renewcommand*{\P}{\mathbb{P}}

\newcommand{\eps}{\varepsilon}
\newcommand{\rr}{\mathbb{R}}

\newcommand{\nn}{\mathbb{N}}

\newcommand{\A}{\mathcal{A}}

\newcommand{\M}{\mathcal{M}}
\newcommand{\I}{\mathcal{I}}
\newcommand{\U}{\mathcal{U}}

\newcommand{\lsc}{\textup{lsc}}
\newcommand{\cl}{\textup{cl}}
\newcommand{\co}{\textup{co}}

\newcommand{\one}{\mathds{1}}

\newcommand{\defined}{\mathrel{\mathop:}}

\DeclareMathOperator*{\argmin}{arg\,min}

\theoremstyle{plain}
\newtheorem{theo}{Theorem}[section]
\newtheorem{lemma}[theo]{Lemma}
\newtheorem{propo}[theo]{Proposition}

\theoremstyle{definition}

\newtheorem{defi}[theo]{Definition}

\newtheorem{remark}[theo]{Remark}
\newtheorem*{standing_assumption}{Standing Assumption (A)}
\newtheorem*{assumption2}{Assumption (B)}

\numberwithin{equation}{section}

\cleardoublepage
\title{On first passage time problems of Brownian motion --\\ The inverse method of images {revisited}}

\author{Sören Christensen\thanks{Kiel University, \emph{Email:} christensen\MVAt math.uni-kiel.de} 
	\and 
	Oskar Hallmann
	\and Maike Klein\thanks{Kiel University, \emph{Email:} maike.klein\MVAt math.uni-kiel.de.
}
}

\begin{document}

	\maketitle
	\begin{abstract}

	Let $W$ be a standard Brownian motion with $W_0 = 0$ and let $b\colon [0,\infty) \to \rr$ be a continuous function with $b(0) > 0$. In this article, we look at the classical First Passage Time (FPT) problem, i.e., the question of {determining} the distribution of $\tau \coloneqq \inf \{ t\in [0,\infty)\colon  W_t \geq b(t) \}.$ More specifically, we revisit the method of images, which we feel has received less attention than it deserves. The main observation of this approach is that the FPT problem is fully solved if a measure $\mu$ exists such that 
		\begin{align*}
			\int_{(0,\infty)} \exp\left(-\frac{\theta^2}{2t}+\frac{\theta b(t)}{t}\right)\mu(d\theta)=1, \qquad t\in(0,\infty).
		\end{align*}
The goal of this article is to lay the foundation for answering the still open question of the existence and characterisation of such a measure $\mu$ for a given curve $b$. We present a new duality approach that allows us to give sufficient conditions for the existence. Moreover, we introduce a very efficient algorithm for approximating the representing measure $\mu$ and provide a rigorous theoretical foundation.
	\end{abstract}

	\begin{center}\footnotesize
		\begin{tabular}{l@{ : }p{9.5cm}}
			{\itshape 2020 MSC} & 90C05, 60G40, 60J65.
			\\
			{\itshape Keywords} & first passage time, inverse method of images, linear programming, Brownian motion, boundary hitting.
		\end{tabular}
	\end{center}

	\section{Introduction}
	
	One of the most intuitive questions to ask is when a Brownian motion crosses a given boundary. This problem is widely known as the \textit{First Passage Time (FPT) problem}. Let $W$ be a one-dimensional Brownian motion with $W_0 = 0$ and let $b\colon [0,\infty) \to \rr$ be a continuous function with $b(0) > 0$. The first passage time of $b$ (from below) is then defined as
	\begin{align*}
		\tau \coloneqq \inf \{ t \in[0,\infty)\colon W_t \geq b(t) \}.
	\end{align*}
	One goal of the FPT problem is to determine the cumulative distribution function (c.d.f.) $F$ of~$\tau$,  
	i.e., $F(t) = \P(\tau \leq t)$.
At the end of this section we give a survey of the literature on the FPT problem, where it turns out that none of the numerous solution methods answers the question completely satisfactorily.  The \textit{method of images} takes a special position among these methods because it can be used to obtain a large family of curves for which the first passage time distribution is known quite explicitly. The application in this context goes back at least to ideas of \cite{levy1965processus} and has been widely used for questions of sequential statistics \cite{DAN82,robbins1985statistical}. We describe the essential ideas now based on \cite{LER86} as far as necessary for understanding and give more details and an intuition in Section \ref{ChapMethodofImages}. 

Let $\mu$ be a positive, $\sigma$-finite measure on $(0,\infty)$ and $b\colon [0,\infty) \to \R$ be the unique solution  $b(t)=x$ of the equation
\begin{align}
	\int_{(0,\infty)} \exp\left(-\frac{\theta^2}{2t}+\frac{\theta x}{t}\right)\mu(d\theta)=1, \qquad t\in(0,\infty).\label{eq:determine_b}
\end{align}
 Then the distribution function of $\tau$, the first passage time of the Brownian motion $W$ to $b$, is given as a function of $\mu$ and $b$, more precisely, 
\begin{align*}
		\P(\tau\leq t)	&= 1-\Phi \left( \frac{b(t)}{\sqrt{t}} \right) + \int_{(0,\infty)}  \Phi \left( \frac{b(t)-\theta}{\sqrt{t}} \right) \mu(d \theta).
\end{align*}
Therefore, the method of images can be seen as a method to generate certain curves with explicitly given first passage time distribution which is why the method has been celebrated. However, explicit solutions remain scarce apart from the linear boundary, cf.\ \cite[pp.27 ff.]{LER86} for more explicit examples.

In practice, it seems more reasonable to ask the inverse question: Given a boundary $b$, how must $\mu$ be chosen to obtain the distribution of $\tau$ using the method of images, i.e., such that $b$ is the solution to \eqref{eq:determine_b}. This is called the \textit{inverse method of images}.
This question has already been formulated in \cite{LER86}  and was more recently raised again in \cite{KAH08}.
	 It should actually be divided into two parts:
\begin{itemize}
	\item Is it possible to identify classes of curves $b$ that ensure the existence of a representing measure~$\mu$ in the sense of the method of images, i.e., $\mu$ satisfies $\int_{(0, \infty)} \exp\left(-\frac{\theta^2}{2t}+\frac{\theta b(t)}{t}\right)\mu(d\theta)=1$ for all $t\in(0,\infty)$? To the best of our knowledge, there has been no significant progress on this question. 
	\item Given that a curve is representable, how can we approximate the measure $\mu$ numerically in an efficient and in a theoretically sound way? Initial approaches to this, albeit without e.g.\ a guarantee of convergence, already exist in the literature. We discuss this in more detail in Section \ref{ChapCompMeth}. 
\end{itemize}
In summary, it can be said that there are only very partial results on the inverse method of images. This may explain why the method has received less attention in recent years than the authors of this article believe it deserves. This paper aims to change this and thus to stimulate new activity in this area. Our main contributions are as follows:
\begin{itemize}
	\item We formulate and investigate the representability of curves by studying two dual infinite dimensional linear programs. This opens up a new perspective on the problem and provides a basis for formulating sufficient conditions for a curve to be representable by a measure $\mu$. 
	\item In addition, the adaptation of infinite linear programming algorithms provides an extremely efficient and accurate way to solve the FPT problem numerically. This method is further underpinned by precise convergence results. 
\end{itemize}

	\subsection{Structure of the paper}

The paper is organised as follows. 
In Section \ref{ChapMethodofImages} we introduce the method of images and the inverse method of images in more detail. 
Section \ref{ChapLinProg} deals with a new infinite  dimensional linear programming approach to the inverse method of images. We give two different linear programs, formulate the corresponding dual programs and prove strong duality results for both set-ups (Sections \ref{SecSD1} and \ref{SecSD2}). 

This forms the basis for our key results. First, we obtain sufficient conditions for representability of a given concave, analytic boundary $b$, see Section \ref{SecExisRepMeasure}.
In Section \ref{ChapCompMeth} we provide convergence results for our linear programs, where we only discretise the space axis which allows to present a new algorithm and  give error bounds for the approximation. In Section \ref{SecNumStuRep} we present the numerical results for our algorithm for four boundaries. Finally, we shortly  comment on the inverse method of images for two-sided boundaries in Section \ref{ChapTwoSidedBoundaries}.

\subsection{Literature on the FPT problem}
	The FPT problem has been studied extensively and many applications have been developed: In statistics, FPT problems arise for problems in testing, cf.\ \cite{RS73, FER82B} or the surveys \cite{SIE86, LAI01}. In finance, the problem emerges in the valuation of barrier options, cf.\ \cite{KI92,GY96,RS97}, as well as in default models, cf.\ \cite{CGJ08,HW01}, and in the evaluation of credit risks, cf.\ \cite{CHE06}. For an overview of the FPT problem and some applications in the field of physics, including a connection to electrostatics, see e.g.\ \cite{RED22} or the monograph \cite{RED01}.

	The FPT problem can be traced back to the thesis ``Théorie de la spéculation'' by Bachelier~\cite{BAC00}, where the problem is first formulated, but only for constant $b$. Other early formulations of the problem include \cite{SCH15,KHI33}. 
	
	Even though the problem has been studied for a long time, closed form solutions for the distribution function $F$ of $\tau$ given the function $b$ are rare. For linear boundaries $b(t) = a+m\:\!t$ with  $a>0$, $m \in \rr$ the well-known Bachelier-Lévy formula gives the density of $\tau$ as
	\begin{align*}
		\frac{a}{t^{3/2}} \,\varphi\left( \frac{a+m\:\!t}{\sqrt{t}} \right), \qquad t\in(0,\infty),
	\end{align*}
	where $\varphi$ is the density of a standard normal distribution, see \cite[Chapter III, §3]{RY99} for the case $m=0$ or \cite[Example 1, p.27]{LER86} for the more general case $m \in \rr$.
	
	Partial results are available for other kind of boundaries: 
	For the square-root boundary, asymptotic results for the distribution function $F$ of $\tau$ are obtained in \cite{BRE67},  a rather lengthy but explicit formula for the c.d.f.\ $F$ associated to square-root boundaries is derived in \cite{RSS84}  and  an infinite power series is presented in \cite{NFK99}. 	
	The quadratic boundary was investigated among others in \cite{SAL88, GRO89}, where in both papers a formula for the distribution function of $\tau$ is derived which depends on Airy functions.
	
	Another approach to the FPT problem is by integral equations (typically Volterra or Fredholm type integral equations) connecting the boundary $b$ and the distribution function $F$ of $\tau$.
	While there have been many approaches involving Volterra type integral equations for many years, Peskir~\cite{PES02} presents a unifying approach via an integral equation -- the so-called ``master equation'' --  to derive these equations.
	These integral equations are in most cases difficult to solve analytically but there are numerical approaches, see e.g.\ \cite{DUR71,DUR85, DUR92, SMI72, PP74,PES02,NAR01}.
		There has also been use of Fredholm type integral equations but more rarely, see \cite{SHE67,NOV81,DAN00,JKV09,MR4635691}.
			
	In recent years, numerical approaches have rather concentrated on Monte Carlo methods, see e.g.\ \cite{WP97,PW01,BN05,Poe12,JW17}.		 Their most important drawbacks are the extensive computation time and the problem of undetected crossings in between discretisation steps. 
	
	\section{The method of images and the inverse method of images}\label{ChapMethodofImages}

\begin{figure}[h!]
	\begin{center}
\includegraphics[width=\textwidth]{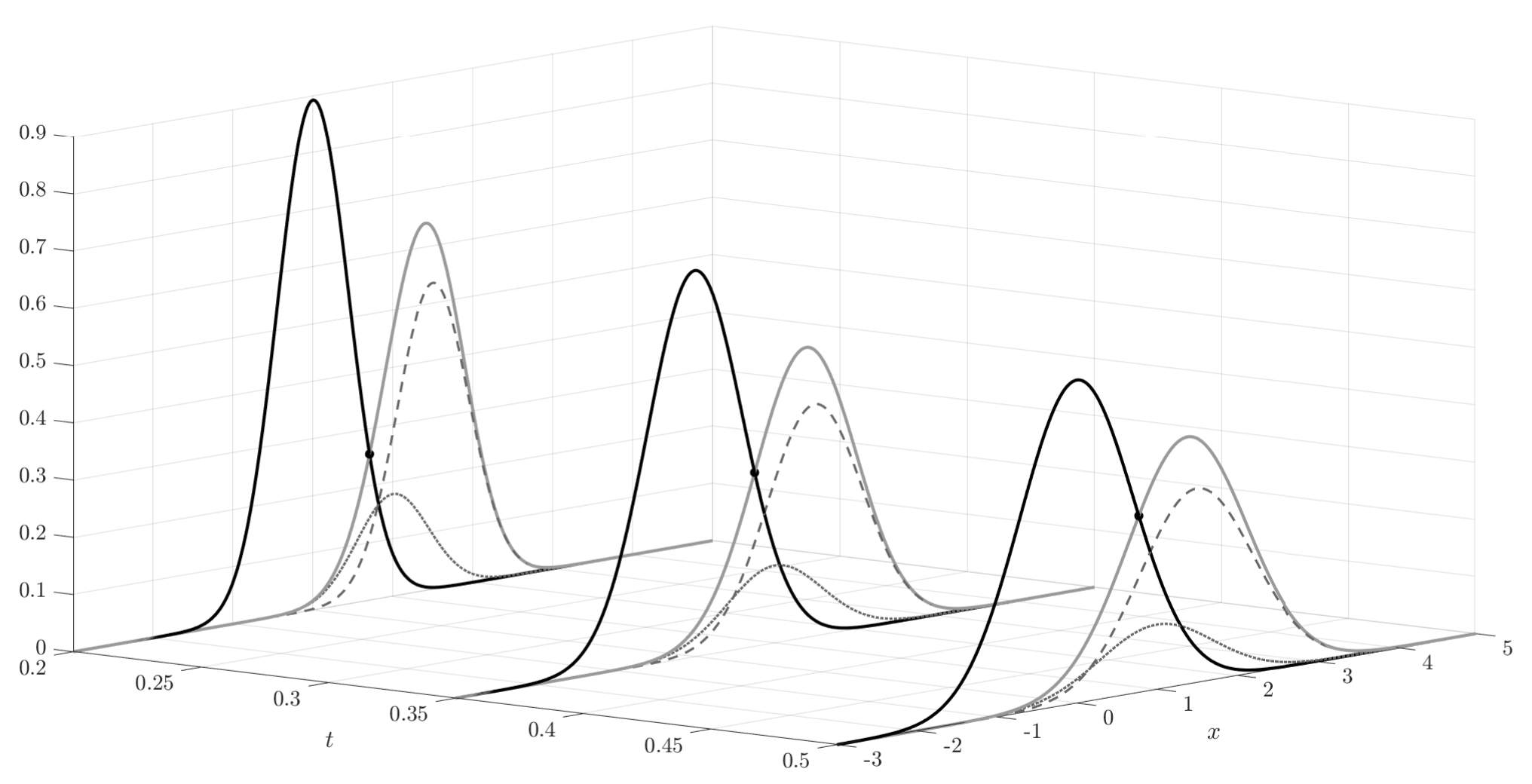}
\includegraphics[width=0.5\textwidth]{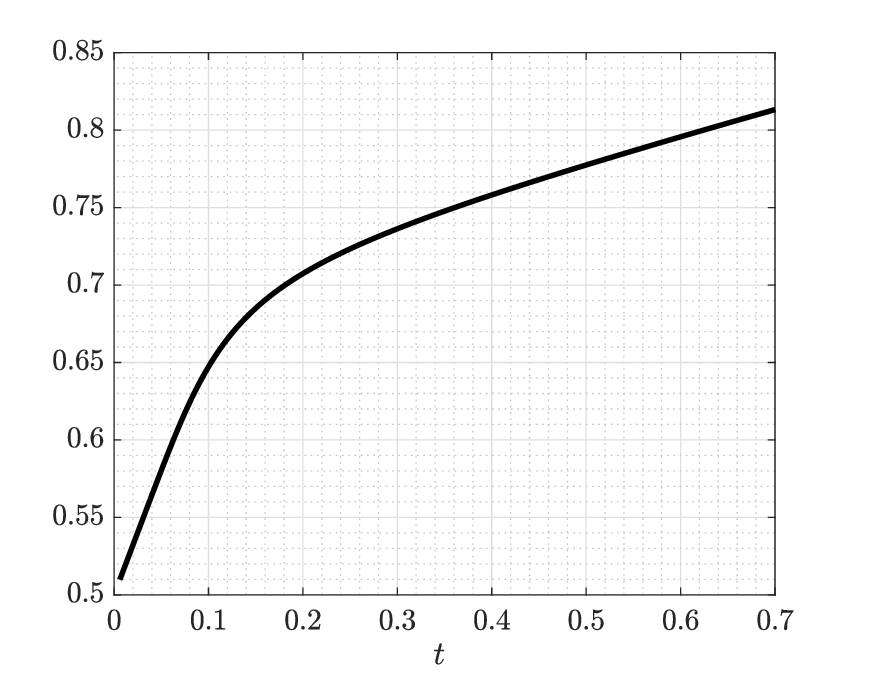}
\end{center}
\caption{The method of images for $\mu=\frac15\delta_1+\frac35\delta_{\frac32}$.\\[0.05cm]
			\emph{Top}: The idea of the method of images for the time points $t\in\{0.2,0.35,0.5\}$: We determine the point $b(t)$ in space, where the density of a Brownian motion at time $t$ started in $0$ (black) coincides with the {$\mu$-weighted densities of two Brownian motions at time $t$ started in $1$ and $\frac32$}, respectively, (solid, gray). The dashed and dotted gray functions are the weighted densities $\frac15\frac{1}{\sqrt{t}}\varphi\left(\frac{\cdot-1}{\sqrt{t}}\right)$ and $\frac35\varphi\left(\frac{\cdot-3/2}{\sqrt{t}}\right)$ which sum up to the $\mu$-weighted density.\\[0.1cm]
	\emph{Bottom}: The boundary $b$ obtained from the method of images for $\mu$.}\label{fig:density}
\end{figure}

	The method of images 
	 is based on the following idea (cf.\ \cite[p.18]{LER86}): consider on the one hand a Brownian motion starting in 0  with unit mass and on the other hand a second Brownian motion with {starting points greater than $0$ given by some mass distribution $\mu$.} 
	 We then observe the ``superposition'' of the distributions of the two Brownian motions, or more precisely the points $(t,x)$ in space-time at which the Brownian motions are staying with the same intensity. 
	 These are the points $(t,x)\in (0,\infty)\times\rr$ such that {the density of the standard Brownian motion at time $t$ evaluated in $x$ coincides with the density of the other Brownian motion at time $t$ weighted according to $\mu$ also evaluated in $x$, see Figure \ref{fig:density}.} The main observation in the method of images is that these points form a boundary $b$ for which the first passage time distribution $F$ can be given explicitly in terms of $\mu$ (see Proposition \ref{prop:meth_images} below).

\newpage
We define
\begin{align*}
	h(t,x) = \frac{1}{\sqrt{t}} \varphi \left( \frac{x}{\sqrt{t}} \right) - \int_{(0,\infty)}  \frac{1}{\sqrt{t}} \varphi \left( \frac{x-\theta}{\sqrt{t}} \right) \mu(d\theta), 
\end{align*}
	where $\varphi$ is the density of the standard normal distribution and
\begin{align*}
		r_\mu(t,x) = r(t,x) = \int_{(0,\infty)}  r_\theta(t,x) \mu(d\theta),
\end{align*}
 with
\begin{align*}
	r_\theta(t,x)=\exp\left(-\frac{\theta^2}{2t}+\frac{\theta x}{t}\right).
\end{align*}
In most cases we will write $r$ instead of $r_\mu$ whenever it is clear which measure $\mu$ we are referring to. Then, we have
\begin{align*}
	h(t,x) = \frac{1}{\sqrt{t}} \varphi \left( \frac{x}{\sqrt{t}} \right) \big(1-r(t,x)\big).
\end{align*}
Now, we consider a given $b\colon [0,\infty) \rightarrow \rr$. The method of images  refers to the case that $\mu$ and $b$ are related in the following way: 
\begin{align}\label{eq:b_sol_to_r}
	r\big(t, b(t)\big)=1,\qquad t\in(0,\infty), 
\end{align}
which is equivalent to $h\big(t,b(t)\big)=0$ for all $t\in(0,\infty)$. It can be easily shown that a unique solution $b$ exists to \eqref{eq:b_sol_to_r} if $r_\mu$ is finite, which is guaranteed if $\mu$ is a positive measure on {$(0,\infty)$
	 such that for all $\varepsilon > 0$
	\begin{align*}
		\int_{(0,\infty)} \varphi(\sqrt{\varepsilon} \theta ) \mu(d\theta) < \infty.
	\end{align*}
}

All $b$ satisfying \eqref{eq:b_sol_to_r} for some $\mu$ have certain properties:
\begin{lemma}\label{LemmaAnalyticConcave}
	Let
	\begin{align*}
		\theta^\ast = \inf\left\{\theta \in(0,\infty)\colon \mu\big((0, \theta]\big) >0\right\} \geq 0.
	\end{align*}
	 Furthermore, let $b\colon [0,\infty)\to\R$ satisfy $r_\mu\big(t, b(t)\big) = 1$ for all $t\in(0,\infty)$. Then,
	\begin{enumerate}[(a)]
		\item\label{item:analytic} $b$ is analytic.
		\item $b$ is concave.
		\item\label{LemmaInitialMass} $\lim_{t \searrow 0} b(t) = \frac{\theta^\ast}{2}.$
	\end{enumerate}
\end{lemma}
\begin{proof}  
	 See the proof of Lemma 1.1 and Lemma 1.2 in  \cite{LER86} and slightly extend the arguments for \ref{item:analytic}.
		\end{proof}

\begin{remark}\label{rem_LemmaAnalyticConcave_on_0_t0}
	If there exists $t_0\in(0,\infty)$ with $r_\mu\big(t, b(t)\big) = 1$ for all $t \in (0,t_0]$ only, then (a) and (b) hold on $(0,t_0)$. Moreover, (c) is also satisfied.
\end{remark}

Now assume that $W$ is a one-dimensional Brownian motion starting in $0$  and let $\tau$ be the first passage time of $W$ to $b$, i.e.,
\begin{align*}
	\tau \coloneqq \inf \big\{ t \in[0,\infty) \colon W_t \geq b(t) \big\}.
\end{align*}
 Here $b$ is again a given curve and we will from now on \emph{not} necessarily assume that $b$ fulfills $r\big(t, b(t)\big)=1$ for a given measure $\mu$. 
  If $\mu$ is chosen in a suitable way, the function $r$ can be used to approximate the distribution of the first passage time of $W$ to $b$.
\begin{propo}\label{prop:meth_images}
	Let $t_0\in(0,\infty)$ and let $b\colon [0,\infty) \to \rr$ be a continuous function such that there exist $\zeta_1\in[0,1)$ and $\zeta_2\geq 0$ with
	\begin{align*}
		\begin{aligned}
		1-\zeta_1\leq \frac{1}{r\big(t,b(t)\big)}\leq 1+\zeta_2, \qquad t\in(0,t_0].
		\end{aligned}
	\end{align*}
	Then, for all $s \in (0, t_0]$ and all $x < b(s)$ it holds
	\begin{align*}
		(1-\zeta_1)\:\!r(s,x)\leq \P(\tau\leq s|W_{s}=x)\leq (1+\zeta_2)\:\!r(s,x).
	\end{align*}
	In particular, if $r\big(t,b(t)\big)=1$ for all $t\in(0,t_0]$, it holds for all $s \in (0, t_0]$ and all $x < b(s)$ that
	\begin{align*}
		r(s,x)=\P(\tau\leq s|W_s=x).
	\end{align*}
\end{propo}
	\begin{proof}
		For the proof slightly extend the arguments in the alternative proof of Theorem 1.1 in \cite[pp.40 f.]{LER86}.
		 For more details we refer to the proof of Proposition 1.3 in \cite{DissOskar}.
\end{proof}

If $b$ satisfies $r\big(t,b(t)\big) = 1$ for all  $t \in (0,T]$ with $T\in (0,\infty)\cup\{\infty\}$ (with the convention that $(0,\infty]=(0,\infty)$), we can recover the distribution function $F(t) = \P(\tau \leq t)$ for $t\in (0,T]$ from $r(t,x) = \P(\tau \leq t \:\!\vert\:\! W_t =x)$ by integration. Indeed,
\begin{align}
	\begin{split}
	F(t) &= \P\big(W_t \geq b(t)\big) + \int_{(-\infty,b(t))} \P(\tau \leq t \:\!\vert\:\! W_t = x) \,p_t(0,x) dx\\
	&= 1-\Phi \left( \frac{b(t)}{\sqrt{t}} \right) + \int_{(0,\infty)}  \Phi \left( \frac{b(t)-\theta}{\sqrt{t}} \right) \mu(d \theta), 
	\end{split}\label{eq:density}
\end{align}
where $p_t(0,x)=\frac{1}{\sqrt{2\pi t}} \exp\big(-x^2/(2t)\big)$ is the transition kernel of a Brownian motion from $0$ at time $0$ to $x$ at time $t$.
 In particular, the density $f$ of $\tau$ is given by
\begin{align*}
	f(t) = \frac{1}{2t^{3/2}} \int_{(0,\infty)} \theta\, \varphi \left( \frac{\theta-b(t)}{\sqrt{t}} \right) \mu(d \theta).
\end{align*}

The \emph{method of images} is traditionally applied as follows: One starts with a measure $\mu$ with associated function $r = r_\mu$ and then considers a curve $b$ that is the implicit solution to the equation
\begin{align*}
	r(t,x)=1,
\end{align*}
i.e., $b$ is chosen such that $r\big(t, b(t)\big) = 1$ for all $t \in (0, t_0]$ for some $t_0 >0$. Proposition \ref{prop:meth_images} then yields that, under certain assumptions, $r$ can be used to derive the passage time distribution of the curve~$b$ using \eqref{eq:density}.
\begin{remark}\label{rem:linear}
	Using the method of images one can generate curves with explicit {hitting probabilities}. The easiest examples are linear boundaries $b(t)=a+m\:\!t$ with $a>0, m \in \rr$, which are generated by $\mu=\exp(-2am)\delta_{2a}$, where $\delta_{2a}$ denotes the Dirac measure in $2a$,  cf.\ \cite[Example 1, p.27]{LER86}. 
\end{remark}

As discussed in the introduction, we are interested in the \emph{inverse method of images}: given a curve $b$, does there exist a measure $\mu$ such that the method of images applied to $\mu$ yields $b$? 
We approach this question in the following section by methods of linear program.

In order to do so, we start the following definition.
\begin{defi}
	A function $b\colon [0,\infty) \to \rr$ is called \emph{representable (in the sense of the method of images) up to $T$} for some $T\in(0,\infty)\cup\{\infty\}$ if there exists a  measure $\mu$
	 such that $b$ is the solution of
	\begin{align*}
		\int_{(0,\infty)}  r_\theta(t,x) \mu(d\theta) = r(t,x) = 1, \qquad t\in(0,T]
	\end{align*}
	i.e., $r\big(t, b(t)\big) = 1$ for all $t\in(0,T]$. In this case, we call $\mu$ the \emph{representing measure} of $b$ up to $T$ and say $\mu$ \emph{represents} $b$ on $[0,T]$. 
\end{defi}
\begin{remark}
			Here we again use the convention that $(0,\infty]=(0,\infty)$.
\end{remark}

Next we  analyse the existence of representing measures.
  
\begin{theo}\label{ThmRepSeq}
	Let $b\colon[0,\infty)\to \rr$  be analytic with $b(0) >0$. Assume that $\mu$ is a positive $\sigma$-finite measure on $(0,\infty)$ 
	with $\mu\big(\big(0,2b(0)\big)\big)=0$ satisfying 
		\begin{align*}
		\int_{(0, \infty)} \exp\left(-\frac{\theta^2}{2 t^\ast} \right) \mu(d\theta) < \infty
	\end{align*}
	for some $t^*\in (0,\infty)$. Furthermore, let $b$  be analytic with $b(0) >0$. Assume that 
		there exists a sequence $(t_n)_{n\in\nn}$ with accumulation point $\tilde{t} \in (0, t^*)$ such that
		\begin{align*}
			\int_{(0,\infty)}  r_\theta \big(t_n, b(t_n)\big) \mu(d\theta) = 1 \qquad \text{for all } n \in \nn.
		\end{align*}
	Then, $b$ is representable by $\mu$ on $[0,	{t^*})$.
\end{theo}
\begin{proof}
	We show that the function 
	\begin{align*}
	t \mapsto r\big(t,b(t)\big)=\int_{(0, \infty)} r_\theta(\big(t,b(t)\big)\mu(d\theta)
\end{align*}
 is analytic on $(0,\hat{t}\,)$ for every $\hat{t}\in (\tilde{t}, t^*)$. Then the
	claim follows directly from the identity theorem for analytic functions, e.g., cf.\ \cite[Theorem III.3.2]{FB09}. 
	
	Let $\hat{t} \in (\tilde{t}, t^\ast)$. For all $t \in (0, \hat{t}\:\!)$ observe that
		\begin{align*}
		r\big(t,b(t)\big) &= \int_{[2b(0), \infty)} \exp\left(-\frac{\theta^2}{2t}+ \frac{\theta b(t)}{t} \right) \mu(d\theta)\\[0.2cm]
		&= \int_{[2b(0), K]} \exp\left(-\frac{\theta^2}{2t}+ \frac{\theta b(t)}{t} \right) \mu(d\theta) + \int_{(K, \infty)} \exp\left(-\frac{\theta^2}{2t}+ \frac{\theta b(t)}{t} \right) \mu(d\theta),
	\end{align*}
	where 
	\begin{align*}
		K \coloneqq \sup_{t \in\left[0, \hat{t}\, \right]} \frac{2t^\ast b(t)}{t^\ast -t}\geq 2b(0) .
	\end{align*}
	Then, $K$ is finite and $\mu$ is $\sigma$-finite by assumption and therefore the first integral is finite. Observe that for  $\theta > K$ we have 
	\begin{align*}
		\exp\left(-\frac{\theta^2}{2t}+ \frac{\theta b(t)}{t} \right) \leq \exp\left(-\frac{\theta^2}{2 t^\ast} \right), \qquad t \in (0, \hat{t}\,),
	\end{align*}
	and so the second integral is finite by assumption. Thus, $t \mapsto r\big(t, b(t)\big)$ is analytic on $(0, \hat{t}\,)$. 
\end{proof}


	\section{A linear programming approach for the inverse method of images}\label{ChapLinProg}
	
	In order to be able to answer the question of the existence of a representing measure, a candidate must be found. Our approach to this is inspired by \cite{CKL22} and consists of identifying $\mu$ as the solution to a linear optimisation problem. More precisely, we investigate two linear programs that give upper and lower bounds for the first passage time distribution of a given boundary $b$. Then, we show duality results for these programs as a basis to give sufficient conditions under which a representing measure $\mu$ for $b$ exists.
	In order to have any chance of representability at all, we require the properties of Lemma~\ref{LemmaAnalyticConcave}. 
	
\begin{standing_assumption}\label{assumptionB}	We fix an analytic, concave boundary $b\colon [0, \infty) \rightarrow \rr$ with $b(0) > 0$ and finite slope at$~0$, i.e., $| b'(0)| < \infty$. 
	\end{standing_assumption}

	For a  locally compact Hausdorff space $X$  let $\M(X)$ be the set of all regular signed measures on $X$ with finite total variation and denote by $\M^+(X)$ the cone of all non-negative measures in $\M(X)$. For $\lambda\in\M(X)$ denote by $\| \lambda \|$ its total variation. For $\lambda\in\M^+(X) $ it holds that $\|\lambda\|=\int_X 1 \:\!d\lambda$. 
	
	 Moreover, we fix a point $(t_0,x_0)$ with $t_0>0$ and $x_0< b(t_0)$. In the light of Proposition \ref{prop:meth_images}, it seems natural to consider the following linear problem
	\begin{align}
			\boxed{	
				\begin{aligned}
			& \text{maximise} & & \int r_\theta(t_0,x_0)\mu(d\theta) \\[0.2cm]
			& \text{subject to} & & \mu \in \M^+\big([2b(0),\infty)\big),\\
			&&& \int r_\theta\big(t,b(t)\big)\mu(d\theta)\leq 1 \qquad \text{for any } t\in (0,t_0],
		\end{aligned} \tag{$D_1$}\label{PROGRAM_D_1}
	}
				\end{align}
		where we approximate the measure $\mu$ representing $b$ ``from below''.
	 
	Alternatively, we could also consider the following linear program
	\begin{align}
		\label{PROGRAM_P_2}
		\boxed{\begin{aligned}
			& \text{minimise} & & \int r_\theta(t_0,x_0)\mu(d\theta) \\[0.2cm]
			& \text{subject to} & & \mu \in \M^+\big([2b(0),\infty)\big),\\
			&&& \int r_\theta\big(t,b(t)\big)\mu(d\theta)\geq 1 \qquad  \text{for any } t\in (0,t_0],
		\end{aligned} \tag{$P_2$}
	}
	\end{align}
	where we approximate the measure $\mu$ representing $b$ ``from above''. We omitted the bounds of integration as we integrate over the whole space $(0,\infty)$ in each case. 
	\begin{remark}
		By Lemma \ref{LemmaAnalyticConcave} \ref{LemmaInitialMass} and Remark \ref{rem_LemmaAnalyticConcave_on_0_t0}  it is enough to consider measures $\mu$ on $[2b(0), \infty)$ for a given boundary $b$.
		\end{remark} 
	
	\begin{remark}
		Observe that both programs are linear but infinite-dimensional.
	\end{remark}
	
	\begin{remark}
	The labelling of the first program as \eqref{PROGRAM_D_1} and the second as \eqref{PROGRAM_P_2} may at first glance be confusing. The notation will make more sense when we consider the associated formal dual problems as well as weak and strong duality.
	\end{remark}
	
	The reader may ask herself which program to prefer over the other, i.e., whether there is a ``natural'' or ``better'' choice for one or the other. While we need both programs in order to find sufficient conditions for a representing measure $\mu$ to exist, the most immediate use is that admissible measures to these programs lead to lower and upper bounds for the probability $\P(\tau\leq t_0\:\!|\:\!W_{t_0}=x_0)$ as the following lemma shows.
\begin{lemma}\label{ThmLinProBounds}
	For each \eqref{PROGRAM_D_1}-admissible $\mu_1$ we find that
	\begin{align*}
		r_{\mu_1}(t_0, x_0) = \int r_\theta(t_0,x_0)\mu_1(d\theta)\leq \P(\tau\leq t_0\:\!|\:\!W_{t_0}=x_0);
	\end{align*}
	and for each \eqref{PROGRAM_P_2}-admissible $\mu_2$ it holds  that
	\begin{align*}
		r_{\mu_2}(t_0, x_0) = \int r_\theta(t_0,x_0)\mu_2(d\theta)\geq \P(\tau\leq t_0\:\!|\:\!W_{t_0}=x_0).
	\end{align*}
	Moreover, if there exists a representing measure $\mu^*$ on $[0,t_0]$, i.e., if
	\begin{align*}
		\int r_\theta\big(t,b(t)\big)\mu^*(d\theta)= 1 \qquad  \text{for any } t\in (0,t_0],
	\end{align*}
	then $\mu^\ast$ is a maximiser in \eqref{PROGRAM_D_1} and a minimiser in \eqref{PROGRAM_P_2}.
\end{lemma}
\begin{proof}
	Let $\mu_1$ be \eqref{PROGRAM_D_1}-admissible. Then,
	\begin{align*}
		r_{\mu_1} \big(t, b(t)\big) = \int r_\theta\big(t,b(t)\big)\mu_1(d\theta)\leq 1 \qquad \text{for any } t\in (0,t_0].
	\end{align*}
	
	Setting $\zeta_1=0$ in Proposition \ref{prop:meth_images}  yields 
	\begin{align*}
		r(t_0, x_0) \leq \P(\tau\leq t|W_{t_0}=x_0).
	\end{align*}
	The second inequality follows the same way.
	
	Now let $\mu^\ast$ be a representing measure up to $t_0$, i.e.,
	\begin{align*}
		\int r_\theta\big(t,b(t)\big)\mu^*(d\theta)= 1 \qquad  \text{for any } t\in (0,t_0].
	\end{align*}
	We conclude again by Proposition \ref{prop:meth_images} that
	\begin{align*}
		r(t_0, x_0) = \int r_\theta (t_0, x_0) \mu^\ast (d \theta) = \P(\tau\leq t|W_{t_0}=x_0)
	\end{align*}
	and so $\mu^\ast$ is a maximiser in \eqref{PROGRAM_D_1} and a minimiser in \eqref{PROGRAM_P_2}.
\end{proof}
\begin{remark}
Since we do not use that $b$ is concave or analytic in the proof of Lemma \ref{ThmLinProBounds}, the linear programs give upper and lower bounds for the first passage time distribution to boundaries $b$ that do not need to be concave or analytic.
\end{remark}
By Lemma \ref{ThmLinProBounds} optimal solutions to the programs \eqref{PROGRAM_D_1} and \eqref{PROGRAM_P_2} are candidates for a measure~$\mu$ representing $b$ on $[0,t_0]$. Thus, we now investigate the existence of optimal solutions to \eqref{PROGRAM_D_1} and \eqref{PROGRAM_P_2}  and whether these solutions fulfil $r\big(t,b(t)\big) =1$ for every $t \in (0,t_0]$. 
For the analysis, we follow the standard procedure for linear problems, namely the study and interpretation of (formal) dual problems. However, due to the infinite-dimensional structure of the problems, the question arises as to the choice of the correct spaces in order to obtain duality results. {It turns out that our choice of 
	spaces is appropriate for our purposes.} We refer to \cite{ROC70} or \cite{ROC74}  for an overview on general linear programming.

For \eqref{PROGRAM_D_1}, the dual program is given by
\begin{align}
	\label{PROGRAM_P_1}
	\boxed{
		\begin{aligned}
			& \text{minimise}   & & \|\lambda\|\\
			& \text{subject to} & &   \lambda \in \M^+\big((0,t_0]\big),\\
			&&& \int r_\theta\big(t,b(t)\big)\lambda(dt)\geq r_\theta(t_0,x_0) \qquad  \text{for any } \theta\in [2b(0),\infty)
		\end{aligned} \tag{$P_1$}
	}
\end{align}
and for \eqref{PROGRAM_P_2}, the dual program is
\begin{align}
	\label{PROGRAM_D_2}
	\boxed{\begin{aligned}
			& \text{maximise}   & & \|\lambda\|\\
			& \text{subject to} & &   \lambda \in \M^+\big((0,t_0]\big),\\
			&&& \int r_\theta\big(t,b(t)\big)\lambda(dt)\leq r_\theta(t_0,x_0) \qquad \text{for any } \theta\in [2b(0),\infty).
		\end{aligned} \tag{$D_2$}
	}
	\end{align}
	\begin{remark}
	The maximising problem is always tagged with a ``$D$'' and the minimising problem gets tagged with a ``$P$''. This is in line with the ``usual'' notation for linear programs where the minimising problem is often regarded as the canonical primal problem and the maximising problem is often called the canonical dual problem. 
\end{remark}

Denote by $d_1$, $p_1$, $d_2$ and $p_2$  the optimal values of \eqref{PROGRAM_D_1}, \eqref{PROGRAM_P_1}, \eqref{PROGRAM_D_2} and \eqref{PROGRAM_P_2}, respectively. Then, we have 
\begin{align}\label{eq:weak_duality}
	d_1 \leq p_1\qquad \text{and} \qquad d_2 \leq p_2,
\end{align}
 respectively, i.e., weak duality holds. Indeed, for any \eqref{PROGRAM_D_1}-admissible~$\mu$ and any \eqref{PROGRAM_P_1}-admissible~$\lambda$ the constraints of the linear program imply
 \begin{align*}
	\|\lambda\|=\int 1 \,\lambda(dt)
	\geq \int\int r_\theta\big(t,b(t)\big)\mu(d\theta)\,\lambda(dt)
	= \int\int r_\theta\big(t,b(t)\big)\lambda(dt)\,\mu(d\theta)
	\geq \int r_\theta(t_0,x_0)\mu(d\theta).
\end{align*}
Hence, $d_1\leq p_1$. A similar argument shows $d_2\leq p_2$.

\bigskip

	\subsection{Strong duality of \eqref{PROGRAM_D_1} and its dual \eqref{PROGRAM_P_1}}\label{SecSD1}
	
	In this section we establish strong duality for the programs \eqref{PROGRAM_D_1} and \eqref{PROGRAM_P_1}, i.e., we show $d_1=p_1$. 
	 
Let $\Omega= (0, t_0]$ and
	\begin{align*}
		\begin{aligned}
			&T\colon \M \big([2b(0), \infty)\big) \rightarrow C(\Omega), \qquad & &T\mu(t) = \int_{[2b(0), \infty)} r_\theta\big(t, b(t)\big) \mu(d \theta), \\[0.2cm]
			&T'\colon \M (\Omega) \rightarrow C_0\big([2b(0), \infty)\big),  & &T'\lambda(\theta) = \int_\Omega r_\theta\big(t, b(t)\big) \lambda (dt),
		\end{aligned}
	\end{align*}
	where $C(I)$ and $C_0(I)$ denote the space of continuous functions on an interval $I$ and the space of  continuous functions on $I$ vanishing at infinity, respectively.

	 In addition, let 
	\begin{align*}
		\langle f, \nu \rangle \coloneqq \int f d\nu
	\end{align*}
	 on $C_0(\Omega) \times \M(\Omega)$ and $C_0\big([2b(0), \infty) \big)\times \M\big([2b(0), \infty)\big)$, respectively.
	 Moreover, define 
	 \begin{align*}
	 	g(\theta) = r_\theta(t_0, x_0).
	 \end{align*}
	  We can now reformulate our programs as follows: 
	\begin{align}
		\boxed{
		\label{PROGRAM_D_1_short}
		\begin{aligned}
			& \text{maximise} & & \langle g, \mu \rangle \\
			& \text{subject to} & & \mu \in \M^+\big([2b(0),\infty)\big),\\
			&&& 1 - T\mu \in C^+(\Omega)
		\end{aligned} \tag{$D_1$}
	}
	\end{align}
	and
	\begin{align}
		\label{PROGRAM_P_1_short}
		\boxed{
		\begin{aligned}
			& \text{minimise} & & \|\lambda\|\\
			& \text{subject to} & & \lambda \in \M^+(\Omega),\\
			&&& T'\lambda - g \in C_0^+\big([2b(0),\infty)\big),
		\end{aligned} \tag{$P_1$}
	}
	\end{align}
	where $C^+$, $C_0^+$ denote the cones of non-negative elements in the spaces $C$ and $C_0$, respectively.

	Investigating strong duality between \eqref{PROGRAM_D_1_short} and \eqref{PROGRAM_P_1_short} in our infinite-dimensional setting, is anything but trivial. One main technical problem here is that the underlying set $\Omega$ is not compact. Therefore, we first consider auxiliary problems:

	\subsubsection{Strong duality of the restricted linear programs}
	We now derive a sequence of measures that are the solutions to problems with weaker constraints than in \eqref{PROGRAM_D_1_short} and \eqref{PROGRAM_P_1_short}. In Section \ref{Sec:strong_duality_D1_P1} we show that these measures converge along a subsequence to a measure $\mu_1 \in \M^+\big([2b(0), \infty)\big)$ and $\lambda_1\in \M^+(\Omega)$ which serve as candidates for optimisers in \eqref{PROGRAM_D_1_short} and \eqref{PROGRAM_P_1_short}, respectively.
	
		To define these restricted linear programs, set for any $\varepsilon \in \Omega$
	\begin{align*}
		\Omega_\varepsilon= [\varepsilon, t_0]
	\end{align*}
	and  define
	\begin{align*}
		T^\prime_\eps\colon \M(\Omega_\eps) \rightarrow C_0\big([2b(0), \infty)\big), \qquad T^\prime_\varepsilon\lambda(\theta) = \int_{\Omega_\varepsilon} r_\theta\big(t, b(t)\big) \lambda(dt).
	\end{align*}
	Note that one difference between $T^\prime_\eps$ and $T'$ is the area of integration. On  $C(\Omega_\varepsilon) \times \M(\Omega_\varepsilon)$ and $C_0\big([2b(0), \infty)\big) \times \M\big([2b(0), \infty)\big)$ we also consider the algebraic pairing, i.e., the bilinear mapping,
	\begin{align*}
		\langle f, \nu \rangle = \int f d\nu,
	\end{align*}
which is real-valued and point separating. Moreover, if $C_0(\Omega_\varepsilon)=C(\Omega_\varepsilon)$ and $C_0\big([2b(0), \infty)\big)$ are endowed with the weak topologies $\sigma(C_0, \M)$ and the spaces $\M(\Omega_\varepsilon)$ and $\M\big([2b(0), \infty)\big)$ with the vague topologies $\sigma(\M, C_0)$ induced by the algebraic pairing, then all four spaces are locally convex Hausdorff spaces. Note that the space $C_0(\Omega_\varepsilon)$ is the continuous dual of $\M(\Omega_\varepsilon)$ and vice versa. The same holds for $C_0\big([2b(0), \infty)\big)$ and $\M\big([2b(0), \infty)\big)$. 

	Fubini's theorem implies that 
	\begin{align*}
		\langle T \mu, \lambda \rangle = \langle \mu, T^\prime _\eps\lambda \rangle
	\end{align*}
	with the slight abuse of notation $\langle \mu, T^\prime _\eps\lambda \rangle:=\langle  T^\prime _\eps\lambda, \mu\rangle$
	for every $\mu \in \M\big([2b(0), \infty)\big)$ and $\lambda \in \M(\Omega_\varepsilon)$. 
	Hence, the operators $T$ and $T^\prime_\eps$ are adjoint operators. Moreover,  $T$ and $T^\prime_\eps$ are $\sigma(\M, C_0)-\sigma(C_0, \M)$-continuous (with respect to the corresponding spaces), see  Lemma 5.17 in \cite{LEN17}.

 Now consider the restricted linear program
	\begin{align}
		\label{PROGRAM_D_1_eps}
		\boxed{
		\begin{aligned}
			& \text{maximise} & & \langle g, \mu \rangle \\
			& \text{subject to} & & \mu \in \M^+\big([2b(0),\infty)\big),\\
			&&& 1 - T\mu \in C^+(\Omega_\varepsilon).
		\end{aligned} \tag{$D_{1, \varepsilon}$}
	}
	\end{align}

To show the existence of a maximiser $\mu_\eps$ in \eqref{PROGRAM_D_1_eps}, we start by considering the formal Lagrange dual problem associated to \eqref{PROGRAM_D_1_eps} which is given by 
	\begin{align}
		\label{PROGRAM_P_1_eps}
		\boxed{
		\begin{aligned}
			& \text{minimise} & & \|\lambda\|\\
			& \text{subject to} & & \lambda \in \M^+(\Omega_\varepsilon),\\
			&&& T^\prime_\eps\lambda - g \in C_0^+\big([2b(0),\infty)\big).
		\end{aligned} \tag{$P_{1, \varepsilon}$}
	}
	\end{align}
	Denote the optimal values of \eqref{PROGRAM_D_1_eps} and \eqref{PROGRAM_P_1_eps} by $d_{1, \varepsilon}$ and $p_{1, \varepsilon}$, respectively. Then similar to \eqref{eq:weak_duality} we see that 
	\begin{align*}
		0 \leq d_{1, \varepsilon} \leq p_{1, \varepsilon},
	\end{align*}
	 i.e., weak duality holds. 
	Now we prove the existence of an optimiser in \eqref{PROGRAM_D_1_eps} and \eqref{PROGRAM_P_1_eps}.
	\begin{lemma}\label{Lemma_exis_1_eps}
		There exist a \eqref{PROGRAM_D_1_eps}-admissible $\mu_{1, \varepsilon}$ such that $d_{1, \varepsilon} = \langle g, \mu_{1, \varepsilon} \rangle$ and a \eqref{PROGRAM_P_1_eps}-admissible $\lambda_{1, \varepsilon}$ such that $p_{1, \varepsilon} = \| \lambda_{1, \varepsilon} \|$, i.e., $\mu_{1,\varepsilon}$ and $\lambda_{1,\varepsilon}$ are optimal in \eqref{PROGRAM_D_1_eps} and \eqref{PROGRAM_P_1_eps}, respectively.
		\end{lemma}
	\begin{proof}
	Let $\psi\colon[2b(0),\infty)\to \R,$ $\theta\mapsto\exp\left(-\theta \frac{b(t_0)-x_0}{t_0}\right)$ and observe that $\psi\in C_0\big([2b(0),\infty)\big)$ since  $x_0< b(t_0)$. 
	Consider the following modified version of \eqref{PROGRAM_D_1_eps}: 
		\begin{align}
			\label{PROGRAM_D_1_eps_mod}
			\boxed{\begin{aligned}
				& \text{maximise} & & \langle \psi, \mu\rangle \\
				& \text{subject to} & & \mu \in \M^+\big([2b(0),\infty)\big),\\
				&&& 1- T_\mathrm{mod}\, \mu \in C^+(\Omega_\varepsilon),
			\end{aligned} \tag{$D_{1, \varepsilon, \mathrm{mod}}$}
		}
		\end{align}
		where  
		\begin{align*}
			T_\mathrm{mod}\colon \M\big([2b(0),\infty)\big)\to C(\Omega_\varepsilon), \qquad T_\mathrm{mod}\,\mu(t) = \int_{[2b(0), \infty)} \frac{r_\theta\big(t, b(t)\big)}{r_\theta(t_0, b(t_0))}\, \mu(d\theta).
		\end{align*}
		Note that $T_\mathrm{mod}$ is  continuous 		and that \eqref{PROGRAM_D_1_eps_mod} is equivalent to \eqref{PROGRAM_D_1_eps}. Moreover, for any \eqref{PROGRAM_D_1_eps_mod}-admissible $\mu$ we have
		\begin{align}
			1 \geq T_\mathrm{mod}\, \mu (t_0)&=\int_{[2b(0), \infty)}  \frac{r_\theta\big(t_0, b(t_0)\big)}{r_\theta(t_0, b(t_0))}\,\mu(d\theta)
		= \| \mu \|.
		\label{eq:D_1_eps_mod_admissible_norm_bounded_by_1}
		\end{align}
		Hence, every \eqref{PROGRAM_D_1_eps_mod}-admissible $\mu$ is contained in $B_{\M([2b(0), \infty))	}(1)$, the $\sigma(\M, C_0)$-closed ball of radius $1$ around $0$ in $\M\big([2b(0),\infty)\big)$.
		
		Instead of solving the maximisation problem \eqref{PROGRAM_D_1_eps_mod} we can equivalently maximise 
		the mapping  $\mu \mapsto \langle \psi, \mu \rangle$ over the set
		\begin{align*}
			C_d^\varepsilon \defined = T_\mathrm{mod}^{-1}\left( 1 - C^+(\Omega_\varepsilon) \right) \cap \M^+\big([2b(0),\infty)\big) \cap B_{\M([2b(0),\infty))}(1).
		\end{align*}
					Here we intersect with the set  $ B_{{\M([2b(0),\infty))}}(1)$ which can be shown to be compact (see below) and allows us to conclude that $C_d^\varepsilon$ itself is compact as a closed subset of a compact set. 
						 Indeed,  first observe that $1 - C^+(\Omega_\varepsilon)$ is closed as it is homeomorphic to the $\sigma(C_0, \M)$-closed cone
		\begin{align*}
			C^+(\Omega_\varepsilon) = 	C_0^+(\Omega_\varepsilon)=\bigcap_{\lambda \in \M^+(\Omega_\varepsilon)} \big\{ f \in C(\Omega_\varepsilon) \colon \langle f, \lambda \rangle \geq 0 \big\}.
		\end{align*}
		The continuity of $T_\mathrm{mod}$ implies that $T_\mathrm{mod}^{-1}\big(1-C^+(\Omega_\varepsilon)\big)$ is $\sigma(\M, C_0)$-closed as well. Rewriting
		\begin{align*}
			\M^+\big([2b(0),\infty)\big) = \bigcap_{f \in C_0^+([2b(0),\infty))} \big\{ \mu \in \M\big([2b(0),\infty)\big) \colon \langle f, \mu \rangle \geq 0 \big\}
		\end{align*}
		yields that $\M^+\big([2b(0),\infty)\big)$ is $\sigma(\M, C_0)$-closed, too. Finally, we find that $B_{\M([2b(0),\infty))}(1)$ is $\sigma(\M, C_0)$-compact due to the Alaoglu-Bourbaki theorem, cf.\ \cite[Theorem 23.5]{MV92}. Hence, $C_d^\eps$ is compact as a closed subset of a compact set. 
		
		Since $\psi\in C_0\big([2b(0),\infty)\big)$ the mapping $\mu \mapsto \langle \psi,\mu \rangle$ is continuous with respect to the topology $\sigma(\M, C_0)$  and thus attains its maximal value $d_{1, \varepsilon}$  on $C_d^\varepsilon$ at some measure $\mu_{\varepsilon, \mathrm{mod}} \in C_d^\varepsilon$. Finally, the optimal value $d_{1, \varepsilon}$ of \eqref{PROGRAM_D_1_eps} is attained at $\mu_{1, \varepsilon}$, where
		\begin{align*}
			\frac{d \mu_{1, \varepsilon}}{d \mu_{\varepsilon, \mathrm{mod}}}(\cdot)   =\frac{1}{r_\cdot\big(t_0,b(t_0)\big)}.
		\end{align*}
		
		Now, turning our attention to \eqref{PROGRAM_P_1_eps}, we define $\tilde{\lambda} = \delta_{t_0}$. $\tilde{\lambda}$ is admissible in \eqref{PROGRAM_P_1_eps} for all $\varepsilon\in[0,t_0]$ since it holds that  
		\begin{align*}
			T^\prime_\eps \tilde{\lambda}(\theta) = \exp\left( -\frac{\theta^2}{2t_0} + \frac{\theta b(t_0)}{t_0} \right) > \exp\left( -\frac{\theta^2}{2t_0} + \frac{\theta x_0}{t_0} \right) = g(\theta),\qquad \text{for all } \theta \in [2b(0), \infty)
		\end{align*}
		as $b(t_0) > x_0$ by assumption. Thus, $\tilde{\lambda}$ is admissible and $\big\| \tilde{\lambda} \big\| = 1$. Therefore, it suffices to minimise the $\sigma(\M, C)$-continuous mapping $\lambda \mapsto \| \lambda \|$ over the set
		\begin{align*}
			C_p^\varepsilon \defined = (T^\prime_\eps)^{-1}\left( g+ C^+_0\big([2b(0), \infty)\big) \right) \cap \M^+(\Omega_\varepsilon) \cap B_{\M(\Omega_\varepsilon)}(1).
		\end{align*}
		The $\sigma(\M, C)$-compactness of $C_p^\varepsilon$ now follows along the same lines as the compactness of $C_d^\varepsilon$. Thus, we  conclude 
		 that the continuous function $\lambda \mapsto \| \lambda \|$ attains its minimum $p_{1, \varepsilon}$ at some measure $\lambda_{1, \varepsilon} \in C_p^\varepsilon$.
	\end{proof}
	
	To prove strong duality we first introduce the Lagrange function 
	\begin{align*}
		L_1\colon \M(\Omega_\varepsilon) \times \M\big([2b(0), \infty)\big) \rightarrow[-\infty, \infty]
	\end{align*}
	 associated with the \eqref{PROGRAM_P_1_eps}-\eqref{PROGRAM_D_1_eps}-duality  defined as
	\begin{align*}
		L_1(\lambda, \mu) \defined = \| \lambda \| + \langle g, \mu \rangle - \langle T^\prime_\eps \lambda, \mu \rangle + \I_{\M^+(\Omega_\varepsilon)}(\lambda) - \I_{\M^+([2b(0), \infty))}(\mu),
	\end{align*}
	where we set $L_1(\lambda, \mu)=-\infty$ for $(\lambda,\mu)\in \M^-(\Omega_\eps)\times\M^-\big([2b(0),\infty)\big)$ with $\M^-(I)=\M(I)\backslash\M^+(I)$ for an interval $I$ and where
	\begin{align*}
		\I_M(x) \defined =
		\begin{cases}
			0, \qquad &\text{if } x \in M,\\
			\infty, &\text{if } x \notin M,
		\end{cases}
	\end{align*}
	for any set $M$. 
	We will use the following simplifications later on:
		\begin{align}
				\begin{split}
		\inf_{\lambda \in \M(\Omega_\varepsilon)} \sup_{\mu \in \M([2b(0), \infty))} L_1(\lambda, \mu)&= \inf_{\lambda \in \M^+(\Omega_\varepsilon)} \sup_{\mu \in \M^+([2b(0), \infty))} \left( \| \lambda \| + \langle g-T^\prime_\eps \lambda, \mu \rangle \right) \\
		&= \inf_{\lambda \in \M^+(\Omega_\varepsilon)} \left( \| \lambda \| + \sup_{\mu \in \M^+([2b(0), \infty))} \langle g - T^\prime_\eps \lambda, \mu \rangle \right)\\
		&= \inf_{\substack{\lambda \in \M^+(\Omega_\varepsilon)\\ T^\prime_\eps \lambda \geq g}} \| \lambda \| = p_{1, \varepsilon},
		\end{split}\label{ExpCalcPrimal}\\
	\begin{split}	
		\sup_{\mu \in \M([2b(0), \infty))} \inf_{\lambda \in \M(\Omega_\varepsilon)} L_1(\lambda, \mu) &= \sup_{\mu \in \M^+\big([2b(0), \infty)\big)} \inf_{\lambda \in \M^+(\Omega_\varepsilon)} \left( \| \lambda \| + \langle g-T^\prime_\eps \lambda, \mu \rangle \right) \\
		&= \sup_{\mu \in \M^+([2b(0), \infty))} \left( \langle g, \mu \rangle + \inf_{\lambda \in \M^+(\Omega_\varepsilon)} \langle \lambda, 1- T\mu \rangle \right)\\
		&= \sup_{\substack{\mu \in \M^+\big([2b(0), \infty)\big)\\ T\mu \leq 1}} \langle g, \mu \rangle = d_{1, \varepsilon}.
		\end{split}	\label{ExpCalcDual}
	\end{align}

	Moreover, define the \emph{dual value function} $v_1$ by 
	\begin{align*}
		v_1(f) \defined = \inf_{\mu \in \M([2b(0), \infty))} L^\ast_{1,\mu} (f), \qquad f\in  C(\Omega_\varepsilon), 
	\end{align*}
	where $L^\ast_{1,\mu}$ with
	\begin{align*} L^\ast_{1,\mu}(f)=\sup_{\lambda\in\M(\Omega_\varepsilon)}\left\{\langle f,\lambda\rangle - L_{1,\mu}(\lambda)\right\}
		\end{align*} 
		is the convex conjugate of the mapping $L_{1,\mu}=L_1(\cdot,\mu)$.

	\begin{lemma}\label{Lemma_biconjugate}
		The dual value function $v_1$ is convex and we have $v_1(0)=-d_{1, \varepsilon}  $ and $v_1^{\ast\ast}(0) = -p_{1, \varepsilon}$, where $v_1^{**}=(v_1^*)^*$ is the convex  biconjugate of $v_1$. In particular, it holds $-d_{1, \varepsilon} \geq -p_{1, \varepsilon}$.
	\end{lemma}
	\begin{proof}
		We find with Lemma 5.25 in \cite{LEN17} and our calculations in \eqref{ExpCalcDual} that
		\begin{align*}
			v_1^{\ast\ast}(0) \leq v_1(0) &= \inf_{\mu \in \M([2b(0), \infty))} L^\ast_{1, \mu}(0)\\
			&= \inf_{\mu \in \M([2b(0), \infty))} \sup_{\lambda \in \M(\Omega_\varepsilon)} \left( \langle 0, \lambda \rangle - L_{1, \mu}(\lambda) \right)\\
			&= - \sup_{\mu \in \M([2b(0), \infty))} \inf_{\lambda \in \M(\Omega_\varepsilon)} L_1(\lambda, \mu) = -d_{1, \varepsilon}.
		\end{align*}
		The conjugate $v_1^\ast\colon \M(\Omega_\varepsilon) \rightarrow [-\infty, \infty]$ of $v_1$ is given by
		\begin{align*}
			v_1^\ast(\lambda) &= \sup_{f \in C(\Omega_\varepsilon)} \left( \langle f, \lambda \rangle - v_1(f) \right)\\
			&= \sup_{\mu \in \M([2b(0), \infty))} \sup_{f \in C(\Omega_\varepsilon)} \left( \langle f, \lambda \rangle - L_{1,\mu}^\ast(f) \right)\\
			&= \sup_{\mu \in \M([2b(0), \infty))} L_{1,\mu}^{\ast\ast}(\lambda)\\
			&= \sup_{\mu \in \M([2b(0), \infty))} L_{1,\mu} (\lambda),
		\end{align*}
		where for the last equality we used the Fenchel-Moreau theorem (see e.g.\ Theorem 5 in \cite{ROC74}) which is applicable since for all $\mu \in \M\big([2b(0), \infty)\big)$ the mapping $\M(\Omega_\varepsilon) \ni \lambda \mapsto L_{1, \mu}(\lambda)$ is closed, i.e., either lower semi-continuous and  $L_{1,\mu}(\lambda)>-\infty$ for all $\lambda$ or $L_{1,\mu}\equiv- \infty$, and convex. 
		
		Moreover, we obtain for the biconjugate	$v_1^{\ast\ast}$ of $ v_1$ 
		\begin{align}\label{ExpCalcv1}
			v_1^{\ast\ast}(f) &= \sup_{\lambda \in \M(\Omega_\varepsilon)} \left( \langle f, \lambda \rangle - v_1^{\ast} (\lambda) \right) \nonumber \\
			&= \sup_{\lambda \in \M(\Omega_\varepsilon)} \inf_{\mu \in \M([2b(0), \infty))} \left( \langle f, \lambda \rangle - L_1(\lambda, \mu) \right).
		\end{align}
		The calculation in \eqref{ExpCalcPrimal}  yields
		\begin{align*}
			v_1^{\ast\ast}(0) &= \sup_{\lambda \in \M(\Omega_\varepsilon)} \inf_{\mu \in \M([2b(0), \infty))} \left( \langle 0, \lambda \rangle - L_1(\lambda, \mu) \right)\\
			&= - \inf_{\lambda \in \M(\Omega_\varepsilon)} \sup_{\mu \in \M([2b(0), \infty))} L_1(\lambda, \mu) = -p_{1, \varepsilon}.
		\end{align*}

		Note that  $v_1(f)>-\infty$ for all $f\in C(\Omega_\varepsilon)$. Indeed, assume by contradiction that there exists $f \in C(\Omega_\varepsilon)$ such that $v_1(f) = - \infty$. Then also $v_1^{\ast\ast}(f) = -\infty$ since $v_1^{\ast\ast} \leq v_1$. Then \eqref{ExpCalcv1} implies that $\sup_\mu L_1(\lambda, \mu) = \infty$ for any $\lambda \in \M(\Omega_\varepsilon)$ and thus $p_{1, \varepsilon} = \infty$. But by Lemma \ref{Lemma_exis_1_eps}  there exists a \eqref{PROGRAM_P_1_eps}-admissible solution and therefore,  $v_1(f)>-\infty$ for all $f\in C(\Omega_\varepsilon)$.
		
		Finally, we show that $v_1$ is convex. To this end, let $\alpha \in (0,1)$ and $f_1, f_2 \in C(\Omega_\varepsilon)$. Note that $L_1$ is concave in its second component. For any measures $\tilde{\mu}, \hat{\mu} \in \M\big([2b(0), \infty)\big)$ it holds that 
		\begin{align*}
			v_1(\alpha f_1 + (1-\alpha) f_2) &= \inf_{\mu \in \M([2b(0), \infty))} L^\ast_{1,\mu}(\alpha f_1 + (1-\alpha) f_2)\\
			&= \inf_{\mu \in \M([2b(0), \infty))} \sup_{\lambda \in \M(\Omega_\varepsilon)} \left( \langle \alpha f_1 + (1-\alpha) f_2, \lambda \rangle - L_1(\lambda, \mu) \right)\\
			&\leq \sup_{\lambda \in \M(\Omega_\varepsilon)} \left( \langle \alpha f_1 + (1-\alpha) f_2, \lambda \rangle - L_1(\lambda, \alpha \tilde{\mu} + (1-\alpha)\hat{\mu})\right)\\
			&\leq \alpha \sup_{\lambda \in \M(\Omega_\varepsilon)} \left( \langle f_1, \lambda \rangle - L_1(\lambda, \tilde{\mu}) \right) + (1-\alpha) \sup_{\lambda \in \M(\Omega_\varepsilon)} \left( \langle f_2, \lambda \rangle - L_1(\lambda, \hat{\mu}) \right)\\
			&=\alpha L_{1,\tilde\mu}(f_1)+(1-\alpha) L_{1,\hat\mu}(f_2).
		\end{align*}
		Taking the infimum over $\tilde{\mu}, \hat{\mu}\in\M\big([2b(0), \infty)\big)$ we obtain the convexity of $v_1$.
			\end{proof}
	\begin{remark}
		From the proof of Lemma \ref{Lemma_biconjugate} we conclude that $v_1\colon C(\Omega_\varepsilon)\to (-\infty,\infty]$. 
	\end{remark}
	\begin{remark}
			Lemma \ref{Lemma_biconjugate} implies that $-d_{1, \varepsilon} \geq -p_{1, \varepsilon}$ which is equivalent to weak duality. Of course, we already showed weak duality above without using the additional structure of the dual value function $v_1$.
			\end{remark}

	Now, we can prove  strong duality.

	\begin{propo}\label{Propo_strong_duality_1}
		There exist a \eqref{PROGRAM_P_1_eps}-admissible $\lambda_{1, \varepsilon}$ and a \eqref{PROGRAM_D_1_eps}-admissible $\mu_{1, \varepsilon}$  such that the optimal values $p_{1, \varepsilon}$ and $d_{1, \varepsilon}$ of \eqref{PROGRAM_P_1_eps} and \eqref{PROGRAM_D_1_eps} are attained by $\lambda_{1, \varepsilon}$ and $\mu_{1, \varepsilon}$, respectively. Strong duality holds  between the programs \eqref{PROGRAM_D_1_eps} and \eqref{PROGRAM_P_1_eps}, i.e., $d_{1,\varepsilon}=p_{1,\varepsilon}$, and we have the following complementary slackness conditions
		\begin{align*}
		T\mu_{1, \varepsilon} &= 1, &&\lambda_{1, \varepsilon}\text{-almost surely on } \Omega_\varepsilon\text{ and}\\	
		T^\prime_\eps \lambda_{1, \varepsilon} &= g, &&\mu_{1, \varepsilon} \text{-almost surely on } [2b(0), \infty). 
					\end{align*}
	\end{propo}

	\begin{proof}
		Lemma \ref{Lemma_exis_1_eps} gives the existence of a \eqref{PROGRAM_P_1_eps}-admissible $\lambda_{1, \varepsilon}$ with $p_1=\|\lambda_{1,\eps}\|$ and of  a \eqref{PROGRAM_D_1_eps}-admissible $\mu_{1, \varepsilon}$ with $d_1=\langle g, \mu_{1,\eps}\rangle$.
		
		By	Lemma \ref{Lemma_biconjugate} strong duality holds if $v_1^{\ast\ast}(0) = v_1(0)$. From the Fenchel-Moreau theorem (see e.g., Theorem 5 in \cite{ROC74}) and the convexity of $v_1$ (see Lemma \ref{Lemma_biconjugate}), we deduce that
		\begin{align*}
			v_1^{\ast\ast}(0) = \cl\big(\co(v_1)\big)(0) = \lsc(v_1)(0),
		\end{align*}
		where $\cl$, $\co$ and $\lsc$ denote the closure, the convex hull and the lower semi-continuous hull of a function, for more details see, e.g., \cite[Chapter 3]{ROC74}. 
		Now, use Lemma 5.25 in \cite{LEN17} to obtain
		\begin{align*}
			v_1^{\ast\ast}(0) = \sup_{O \in \U(0)} \inf_{f \in O } v_1(f),
		\end{align*}
		where $\U(0)$ is the set of all $\sigma(C_0, \M)$-open neighbourhoods of $0$. Hence, we want to show that
		\begin{align*}
			v_1(0) = \sup_{O \in \U(0)} \inf_{f \in O } v_1(f),
		\end{align*}
		which is in particular satisfies if $v_1$ is continuous in $0$ with respect to the topology $\sigma(C_0, \M)$.
		Consider $O \defined = \big\{ f\in C_0(\Omega_\varepsilon)\colon \| f \|_\infty < 1\big\}$, which is a $\sigma(C_0, \M)$-open neighbourhood of $0$.  Then, for any $f \in O$ it holds that
		\begin{align*}
			v_1(f) &= \inf_{\mu \in \M([2b(0), \infty))} L^\ast_{1, \mu}(f)\\
			& = \inf_{\mu \in \M([2b(0), \infty))}\ \sup_{\lambda \in \M(\Omega_\varepsilon)} \big( \langle f, \lambda \rangle - L_{1,\mu}(\lambda) \big)\\
			&= \inf_{\mu \in \M^+([2b(0), \infty))}\ \sup_{\lambda \in \M^+(\Omega_\varepsilon)} \left( \langle f, \lambda \rangle - \| \lambda\| - \langle g, \mu \rangle + \langle T^\prime_\eps \lambda, \mu \rangle \right)\\
			&\leq \sup_{\lambda \in \M^+(\Omega_\varepsilon)} \left( \| f \|_\infty \, \| \lambda \| - \| \lambda \| \right) = 0. 
		\end{align*}
	 		Now, 
		 \cite[Theorem 5.42]{AB06} yields that $v_1$ is continuous in $0$ and therefore 
		\begin{align}
			p_{1, \varepsilon} = -v_1^{\ast\ast}(0) = -v_1(0) = d_{1, \varepsilon}.\label{eq:strong_duality}
			\end{align}
		
		Moreover, the optimisers $\lambda_{1, \varepsilon}$ and $\mu_{1, \varepsilon}$ in \eqref{PROGRAM_P_1_eps} and \eqref{PROGRAM_D_1_eps}, respectively,  fulfil the complementary slackness conditions. Indeed, using \eqref{eq:strong_duality} and that $T$ and $T^\prime_\eps$ are adjoint operators   we find that
		\begin{align*}
			0 &\leq \langle T^\prime_\eps\lambda_{1, \varepsilon} - g, \mu_{1, \varepsilon} \rangle = \langle T^\prime_\eps\lambda_{1, \varepsilon}, \mu_{1, \varepsilon} \rangle - d_{1, \varepsilon}= \langle \lambda_{1, \varepsilon}, T \mu_{1, \varepsilon} \rangle - p_{1, \varepsilon} = \langle \lambda_{1, \varepsilon}, T\mu_{1, \varepsilon} - 1 \rangle \leq 0.
		\end{align*}
		In particular, $T^\prime_\eps \lambda_{1, \varepsilon} = g$ holds $\mu_{1, \varepsilon}$-almost surely on $[2b(0), \infty)$ and $T\mu_{1, \varepsilon} = 1$ holds $\lambda_{1, \varepsilon}$-almost surely on $\Omega_\varepsilon$.
	\end{proof}

	\subsubsection{Strong duality of the unrestricted programs }\label{Sec:strong_duality_D1_P1}
	
	We now return to the original programs \eqref{PROGRAM_D_1_short} and \eqref{PROGRAM_P_1_short}. We consider suitable sequences of solutions
	of the restricted programs \eqref{PROGRAM_D_1_eps} and \eqref{PROGRAM_P_1_eps} for $\varepsilon \searrow 0$ which then allow to prove strong duality and the existence of optimisers {(under some assumptions).} 
	
	We now consider measures $\lambda_{\varepsilon}\in \M(\Omega_\eps)$ as measures on $\Omega=(0, t_0]$ by continuing them as the null measure outside of $\Omega_\varepsilon$.

	Recall that 
	\begin{itemize}
		\item in \eqref{PROGRAM_P_1_eps} we can restrict to admissible $\lambda$ with $\| \lambda \| \leq 1$ and
		\item  for the modified program \eqref{PROGRAM_D_1_eps_mod} any admissible  $\tilde{\mu}$ satisfies $\| \tilde{\mu} \| \leq 1$, see \eqref{eq:D_1_eps_mod_admissible_norm_bounded_by_1}. \\
		Any admissible $\mu$ of the unmodified program \eqref{PROGRAM_D_1_eps} can be recovered from an admissible $\tilde{\mu}$ of the modified program by $\frac{d\mu}{ d\tilde{\mu}}(\cdot) =\frac{1}{r_\cdot(t_0,b(t_0))}$, cf.\ the proof of Lemma \ref{Lemma_exis_1_eps}. 
		\end{itemize}
		In particular, the bounds on the total variation are independent of $\varepsilon$. Observe that the metrisation of the vague topology is possible on the total variation unit balls in both spaces (for example, cf.\ \cite[§31]{BAU01}). From the Alaoglu-Bourbaki theorem, cf.\ \cite[Theorem 23.5]{MV92}, we deduce that these unit balls are vaguely compact. Thus, there exists a sequence  $\varepsilon_n \downarrow 0$ and measures $\tilde{\mu}_1$ and $\lambda_1$ with $\| \tilde{\mu}_1 \| \vee \| \lambda_1 \| \leq 1$ such that $\tilde{\mu}_{1, \varepsilon_n} \rightarrow \tilde{\mu}_1$ and $\lambda_{1, \varepsilon_n} \rightarrow \lambda_1$ vaguely as $n\to\infty$.
	
	Moreover, since for every function $f\colon[2b(0), \infty) \to \infty$ with compact support, the function $f/r_\cdot(t_0,b(t_0))$ also has compact support, we conclude that 
	\begin{align*}
		\mu_{1, \varepsilon_n} \xrightarrow[n \to \infty]{} \mu_1, \qquad \text{ vaguely,}
	\end{align*}
	where
	\begin{align*}
		\frac{d\mu_{1, \varepsilon_n}}{d\tilde{\mu}_{1, \varepsilon_n}}(\cdot) :=\frac{1}{r_\cdot\big(t_0,b(t_0)\big)}, \qquad  \frac{d\mu_1}{d\tilde{\mu}_{1}}(\cdot)= \frac{1}{r_\cdot\big(t_0,b(t_0)\big)}. 
		\end{align*}	
	Thus, we have vague convergence of the optimal measures for the restricted measures to some measures $\mu_1$ and $\lambda_1$. These measures are good candidates for being the optimal measures in \eqref{PROGRAM_D_1_short} and \eqref{PROGRAM_P_1_short}. However, we do not yet know whether $\mu_1$ and $\lambda_1$ are even admissible in the respective programs, much less whether these solutions are optimal or whether strong duality holds.
	
	Let us start by looking at admissibility. We know that $\mu_{1, \varepsilon_n}$ converges vaguely to $\mu_1$ and that the mapping $\theta \mapsto r_\theta\big(t, b(t)\big)$ is continuous on $[2b(0), \infty)$ for any $t \in \Omega$ and vanishes at infinity. Thus, by \cite[Theorem 30.6]{BAU01} we have
	\begin{align*}
		T \mu_1 (t) = \int_{[2b(0), \infty)} r_\theta\big(t, b(t)\big) \mu_1(d\theta) = \lim_{n \rightarrow \infty} \int_{[2b(0), \infty)} r_\theta\big(t, b(t)\big) \mu_{1, \varepsilon_n}(d\theta) \leq 1
	\end{align*}
	and so $\mu_1$ is indeed \eqref{PROGRAM_D_1_short}-admissible.
	
	For $\lambda_1$ to be \eqref{PROGRAM_P_1_short}-admissible, we additionally assume that $\| \lambda_{1, \varepsilon_n} \| \rightarrow \| \lambda_1 \|$ as $n\to\infty$. Since then the Portemanteau theorem (see e.g., \cite[Theorem 13.16]{KLE14}) implies that $\lambda_{1, \varepsilon_n} \rightarrow \lambda_1$ weakly and not only vaguely. The only thing that could go wrong is that $\lambda_{1, \varepsilon_n}$ pushes mass into $0$ for $n \rightarrow \infty$ since the total mass of $\lambda_1$ would then be smaller than the mass of the $\lambda_{1, \varepsilon_n}$. In an application the assumption  $\| \lambda_{1, \varepsilon_n} \| \rightarrow \| \lambda_1 \|$ as $n\to\infty$ can be numerically investigated by computing $\lambda_{1, \varepsilon}$ for smaller and smaller $\varepsilon$.
		
	Moreover, in order to guarantee that $t \mapsto r_\theta\big(t,b(t)\big)$ is bounded on $(0, t_0]$ for all $\theta\in[2b(0),\infty)$, we use the assumption that $b'(0) < \infty$. Indeed, for  $\theta > 2b(0)$ we find that
	\begin{align}\label{r_theta_bounded}
		r_\theta\big(t, b(t)\big) = \exp\left( -\frac{\theta}{2t} \big(\theta-2b(t)\big) \right) \xrightarrow[t\downarrow 0]{}\ 0
	\end{align}
 since $\theta - 2b(t) > 0$ for $t$ sufficiently small by the continuity of $b$. For  $\theta = 2b(0)$ it holds that
	\begin{align*}
		r_{2b(0)}\big(t, b(t)\big) = \exp\left( 2b(0) \left( \frac{b(t)-b(0)}{t} \right) \right) \xrightarrow[t\downarrow 0]{} \ \exp \big( 2b(0)b'(0) \big) < \infty
	\end{align*}
	 since we assumed $0 < b(0) < \infty$ and $b'(0) < \infty$. In both cases we find that $t \mapsto r_\theta \big(t, b(t)\big)$ is bounded. Then, the Portmanteau theorem yields that
	\begin{align*}
		\int_{(0, t_0]} r_\theta\big(t, b(t)\big) \lambda_1(dt) = \lim_{n \rightarrow \infty} \int_{(0, t_0]} r_\theta\big(t, b(t)\big) \lambda_{1, \varepsilon_n}(dt) \geq g(\theta)
	\end{align*}
	and so $\lambda_1$ is \eqref{PROGRAM_P_1_short}-admissible.
	
	Finally,  \cite[Theorem 30.6]{BAU01} implies $\lim_{n \rightarrow \infty} \langle g, \mu_{1, \varepsilon_n} \rangle = \langle g, \mu_1 \rangle$. Moreover, by our assumption we have $\lim_{n \rightarrow \infty} \| \lambda_{1, \varepsilon_n} \| = \| \lambda_1 \|$. Thus,
	\begin{align*}
		p_1 &\leq \| \lambda_1 \| \leq \lim_{n \rightarrow \infty} \| \lambda_{1, \varepsilon_n} \| = \lim_{n \rightarrow \infty} p_{1, \varepsilon_n}= \lim_{n \rightarrow \infty} d_{1, \varepsilon_n} = \lim_{n \rightarrow \infty} \langle g, \mu_{1, \varepsilon_n} \rangle = \langle g, \mu_1 \rangle \leq d_1 \leq p_1, 
	\end{align*}
	where we use that weak duality (cf.\ \eqref{eq:weak_duality}) holds in the last inequality.
		In particular,
	\begin{itemize}
		\item $\lambda_1$ and $\mu_1$ are \eqref{PROGRAM_P_1_short}- and \eqref{PROGRAM_D_1_short}-optimal, respectively,
		\item  $d_1 = p_1$, i.e., strong duality holds. 
	\end{itemize}
		Finally, we again use the adjointness of $T$ and $T'$ as well as  $d_1=p_1$ to obtain
	\begin{align*}
		0 \leq \langle T' \lambda_1 - g, \mu_1 \rangle = \langle T' \lambda_1, \mu_1 \rangle - d_1 = \langle \lambda_1, T \mu_1 \rangle - p_1 = \langle \lambda_1, T\mu_1 - 1 \rangle \leq 0
	\end{align*}
	and therefore we have $T^\prime\lambda_1 = g$ holds $\mu_1$-almost surely on $[2b(0), \infty)$ and $T\mu_1 = 1$ holds $\lambda_1$-almost surely on $\Omega$. We summarise our findings.
	\begin{theo}\label{ThmSD1}		
		There exists a sequence $(\varepsilon_n)_{n \in \nn}$ and measures $\lambda_1 \in \M^+(\Omega)$ and $\mu_1 \in \M^+\big([2b(0), \infty)\big)$ such that $\lambda_{1, \varepsilon_n} \rightarrow \lambda_1$ and $\mu_{1, \varepsilon_n} \rightarrow \mu_1$ vaguely. Moreover, $\mu_1$ is \eqref{PROGRAM_D_1_short}-admissible.
			
		In addition, if $\| \lambda_{1, \varepsilon_n} \| \rightarrow \| \lambda_1 \|$, 
		then $\lambda_1$ is \eqref{PROGRAM_P_1_short}-admissible, the optimal values $d_1$ and $p_1$ in \eqref{PROGRAM_D_1} and \eqref{PROGRAM_P_1} are attained by $\mu_1$ and $\lambda_1$, respectively, the optimal values coincide (i.e., strong duality holds) and the following complementary slackness conditions are satisfied: 
			\begin{align*}
				\int r_\theta\big(t,b(t)\big)\mu_1(d\theta)&= 1  &&\hspace{-2cm}\text{for $\lambda_1$-almost all\ } t\in (0,t_0],\\
				\int r_\theta\big(t,b(t)\big)\lambda_1(dt)&=r_\theta(t_0,x_0)  &&\hspace{-2cm}\text{for $\mu_1$-almost all\ } \theta\in [2b(0),\infty).
			\end{align*}
		\end{theo}

	\subsection{Strong duality of \eqref{PROGRAM_P_2} and its dual \eqref{PROGRAM_D_2}}\label{SecSD2}
	As in Section \ref{SecSD1} we consider the restricted linear program
	
	\begin{align}
		\label{PROGRAM_P_2_eps}
		\boxed{\begin{aligned}
			& \text{minimise} & & \langle g, \mu \rangle \\
			& \text{subject to} & & \mu \in \M^+\big([2b(0),\infty)\big),\\
			&&& T\mu - 1 \in C^+(\Omega_\varepsilon)
		\end{aligned} \tag{$P_{2, \varepsilon}$}
	}
	\end{align}
	and its formal dual
	\begin{align}
		\label{PROGRAM_D_2_eps}
		\boxed{
		\begin{aligned}
			& \text{maximise} & & \|\lambda\|\\
			& \text{subject to} & & \lambda \in \M^+(\Omega_\varepsilon),\\
			&&& g - T^\prime_\eps \lambda \in C^+_0\big([2b(0),\infty)\big).
		\end{aligned} \tag{$D_{2, \varepsilon}$}
	}
	\end{align}
	 Denote the optimal values of \eqref{PROGRAM_P_2_eps} and \eqref{PROGRAM_D_2_eps} by $p_{2, \varepsilon}$ and $d_{2, \varepsilon}$, respectively.
	 Again we have weak duality, i.e.,
	 \begin{align*}
	 	d_{2,\eps}\leq p_{2,\eps}.
	 	\end{align*} 

\subsubsection{Strong duality of the restricted linear programs}
	Now we show that there exist optimal solutions to the restricted primal and dual problem, compare Lemma \ref{Lemma_exis_1_eps}. Moreover, strong duality holds. 

\begin{propo}
	There exist a \eqref{PROGRAM_P_2_eps}-admissible $\mu_{2, \varepsilon}$ and a \eqref{PROGRAM_D_2_eps}-admissible $\lambda_{2, \varepsilon}$  such that the optimal values $p_{2, \varepsilon}$ and $d_{2, \varepsilon}$ of \eqref{PROGRAM_P_2_eps} and \eqref{PROGRAM_D_2_eps} are attained by $\mu_{2, \varepsilon}$ and $\lambda_{2, \varepsilon}$, respectively. Strong duality holds, i.e., $d_{2,\varepsilon}=p_{2,\varepsilon}$, and we have the following complementary slackness conditions
	\begin{align*}
		&T\mu_{2, \varepsilon} = 1, \qquad  \lambda_{2, \varepsilon} \text{-almost surely on } \Omega_\varepsilon,\\
		&T^\prime_\eps \lambda_{2, \varepsilon} = g, \qquad  \mu_{2, \varepsilon} \text{-almost surely on } [2b(0), \infty).
	\end{align*}
\end{propo}

\begin{proof} 
	To show the existence of an optimal solution to \eqref{PROGRAM_P_2_eps} and \eqref{PROGRAM_D_2_eps}, respectively, we proceed similar to the proof of Lemma \ref{Lemma_exis_1_eps}. Note that it is enough to show that there exists a \eqref{PROGRAM_P_2_eps}-admissible $\mu$ and that all \eqref{PROGRAM_D_2_eps}-admissible $\lambda$ are bounded by some constant. 
		
	Consider $\bar{\mu} = c \cdot \delta_{2b(0)}$ with
	\begin{align*}
		c > \exp\left( -2b(0)\, \frac{b(t_0)-b(0)}{t_0} \right).
	\end{align*}
	Then, for any $t \in \Omega_\varepsilon$
	\begin{align*}
		T\bar{\mu}(t) &= c \cdot \exp\left( 2b(0) \frac{b(t)-b(0)}{t_0} \right)
		\geq c \cdot \exp\left( 2b(0) \frac{b(t_0)-b(0)}{t_0} \right) > 1,
	\end{align*}
	by the concavity of $b$. Thus,  $T\bar{\mu} - 1 \in C^+(\Omega_\varepsilon)$ and hence, $\bar\mu$ is \eqref{PROGRAM_P_2_eps}-admissible with $\| \bar \mu \| = c < \infty$.
	
	Observe that for any dual feasible $\lambda$ and for $\theta = 2b(0)$ the constraint yields
	\begin{align*}
		g\big(2b(0)\big)
		&= \int_{\Omega_\varepsilon} \exp\left( 2b(0)  \frac{b(t)-b(0)}{t} \right)  \lambda(dt)\\
		&\geq \int_{\Omega_\varepsilon} \exp\left( 2b(0) \frac{b(t_0)-b(0)}{t_0}  \right) \lambda(dt)
		=\exp\left( 2b(0) \left( \frac{b(t_0)-b(0)}{t_0} \right) \right) \| \lambda \|,
	\end{align*}
	by the concavity of $b$. In particular, for any feasible $\lambda$ we have that
	\begin{align*}
		\| \lambda \| \leq \exp\left( 2b(0) \frac{x_0-b(t_0)}{t_0} \right) = \defined C.
	\end{align*}

	In order to show strong duality in the same way as in the first set-up, we consider modified programs. Strong duality for the modified programs follows with the same arguments as in the first set-up and the results carry over to the unmodified programs \eqref{PROGRAM_P_2_eps} and \eqref{PROGRAM_D_2_eps}. For more details we refer to \cite[Section 2.3]{DissOskar}.

		The modification of the programs are necessary for our lines of proof (cf.\ the proof of Proposition~\ref{Propo_strong_duality_1}): our dual value function for the unmodified programs \eqref{PROGRAM_P_2_eps} and \eqref{PROGRAM_D_2_eps} would have been
		\begin{align*}
			v_{2,\mathrm{unmod}}(f) &= \inf_{\lambda \in \M^+(\Omega_\varepsilon)} \sup_{\mu \in \M^+([2b(0), \infty))} \left( \langle f, \mu \rangle - \langle g, \mu \rangle - \| \lambda\| + \langle T\, \mu, \lambda \rangle \right)\\
			&\leq \sup_{\mu \in \M^+([2b(0), \infty))} \left( \| f \|_\infty \| \mu \| - \langle g, \mu \rangle \right).
		\end{align*}
			We could not have defined a $\sigma(C_0, \M)$-open neighbourhood $O$ of $0$ with  $f\leq g$ for all $f\in O$ which then guarantees that $v_{2,\mathrm{unmod}}(f)$ is bounded on $O$. This necessitated our detour.
	\end{proof}
	
	\subsubsection{Strong duality of the unrestricted programs}
	
	As before, we lift the results about the restricted programs \eqref{PROGRAM_P_2_eps} and \eqref{PROGRAM_D_2_eps} to the unrestricted programs \eqref{PROGRAM_P_2} and \eqref{PROGRAM_D_2}. Recall that $p_2$ and $d_2$ denote the optimal values of the unrestricted programs and that weak duality holds, i.e.,
	\begin{align*}
	0 \leq d_2 \leq p_2.
	\end{align*}

    For any $\varepsilon > 0$ the optimisers $\mu_{2,\eps}$ and $\lambda_{2,\eps}$ in \eqref{PROGRAM_P_2_eps} and \eqref{PROGRAM_D_2_eps} satisfy
   \begin{align*}
   	\| \mu_{2, \varepsilon} \|& \leq c\\
   	 		\| \lambda_{2, \varepsilon} \| &\leq \exp\left( 2b(0) \frac{x_0 - b(t_0)}{t_0} \right) = \defined C
	\end{align*}
	for some $c> \exp\left( -2b(0)\, \frac{b(t_0)-b(0)}{t_0} \right)$. In particular, the bounds $c$ and $C$ are independent of $\eps$, thus, with the same metrisation argument as in the first set-up there exists a sequence $(\varepsilon_n)_{n\in\mathbb{N}}$ and measures $\mu_2\in \M^+\big([2b(0),\infty)\big)$ with $\| \mu_2 \| \leq c$ and $\lambda_2\in \M^+(\Omega)$ with $\| \lambda_2 \| \leq C$ such that 
	\begin{align*}
	\mu_{2, \varepsilon_n}\xrightarrow[n\to\infty]{}\ \mu_2,\quad  \text{vaguely,}\\
	\lambda_{2, \varepsilon_n}\xrightarrow[n\to\infty]{}\ \lambda_2, \quad  \text{vaguely.}
	\end{align*}
	
	Since for any $t \in \Omega$ the mapping $\theta \mapsto r_\theta\big(t, b(t)\big)$ is continuous and vanishes at infinity, it follows from \cite[Theorem 30.6]{BAU01}  that
	\begin{align*}
		T \mu_2(t) = \int_{[2b(0), \infty)} r_\theta\big(t, b(t)\big) \mu_2 (d\theta) = \lim_{n \rightarrow \infty} \int_{[2b(0), \infty)} r_\theta \big(t, b(t)\big) \mu_{2, \varepsilon_n} (d \theta) \geq 1.
	\end{align*}
	Hence, $\mu_2$ is \eqref{PROGRAM_P_2}-admissible. 
	
	To show that $\lambda_2$ is \eqref{PROGRAM_D_2}-admissible, observe that for any $\kappa \in \big(0,\frac{t_0}{4}\big)$ we have $\emptyset \neq \Omega_{2\kappa} \subset \Omega_\kappa \subset \Omega$ and so by Urysohn's lemma (e.g., cf.\ \cite[Theorem 4.2]{LAN93}) there exists  a continuous function $\varphi^\kappa\colon \Omega \rightarrow [0,1]$ such that $\varphi^\kappa(t) = 1$ for all $t \in \Omega_{2 \kappa}$ and $\varphi^\kappa(t) = 0$ for all $t \in \mathrm{cl}(\Omega \setminus \Omega_\kappa)$. Since $t \mapsto r_\theta\big(t, b(t)\big)\varphi^\kappa(t)$ is continuous for any $\theta \in [2b(0), \infty)$, we conclude that for any $\theta \in [2b(0), \infty)$
	\begin{align*}
		\int_{(0, t_0]} r_\theta\big(t, b(t)\big) \lambda_2(dt) &= \lim_{\kappa \downarrow 0} \int_{(0, t_0]} r_\theta\big(t, b(t)\big) \one_{\Omega_{2\kappa}}(t) \lambda_2(dt)\\
		&\leq \lim_{\kappa \downarrow 0} \int_{(0, t_0]} r_\theta\big(t, b(t)\big)\varphi^\kappa(t) \lambda_2(dt)\\
		&= \lim_{\kappa \downarrow 0} \lim_{n \rightarrow \infty} \int_{(0, t_0]} r_\theta\big(t, b(t)\big)\varphi^\kappa(t) \lambda_{2, \varepsilon_n}(dt)\\
		&\leq \limsup_{\kappa \downarrow 0} \limsup_{n \rightarrow \infty} \int_{(0, t_0]} r_\theta\big(t, b(t)\big) \lambda_{2, \varepsilon_n}(dt) \leq g(\theta),
	\end{align*}
	where we use \cite[Theorem 30.6]{BAU01} in the second to last line. Therefore,  $\lambda_2$ is \eqref{PROGRAM_D_2}-admissible.

	Since $g$ is continuous and vanishes at infinity,  \cite[Theorem 30.6]{BAU01} yields $\langle g, \mu_2 \rangle = \lim_{n \rightarrow \infty} \langle g, \mu_{2, \varepsilon_n} \rangle$ and therefore
	\begin{align*}
		d_2 \leq p_2 \leq \langle g, \mu_2 \rangle = \lim_{n \rightarrow \infty} \langle g, \mu_{2, \varepsilon_n} \rangle = \lim_{n \rightarrow \infty} p_{2, \varepsilon_n} = \lim_{n \rightarrow \infty} d_{2, \varepsilon_n} \leq d_2.
	\end{align*}
	The last inequality is true as any \eqref{PROGRAM_D_2_eps}-feasible $\lambda$ is also \eqref{PROGRAM_D_2}-feasible. Thus, we conclude that $\mu_2$ is \eqref{PROGRAM_P_2}-optimal,  $d_2 = p_2$, i.e., strong duality holds and $\lim_{n \rightarrow \infty} d_{2, \varepsilon_n} = d_2$.
	
	To show that $\lambda_2$ is \eqref{PROGRAM_D_2}-optimal (and not only admissible) we assume -- similar to the first set-up -- that $\| \lambda_{2, \varepsilon_n} \| \rightarrow \| \lambda_2 \|$ as $n\to\infty$.  Then, 
	\begin{align*}
		\| \lambda_2 \| = \lim_{n \rightarrow \infty} \| \lambda_{2, \varepsilon_n} \| = \lim_{n \rightarrow \infty} d_{2, \varepsilon_n} = d_2.
	\end{align*}
 Finally, we have
	\begin{align*}
		0 \leq \langle T' \lambda_2 - g, \mu_2 \rangle = \langle T'\lambda_2, \mu_2 \rangle - d_2 = \langle \lambda_2, T \mu_2 \rangle - p_2 = \langle \lambda_2, T\mu_2 - 1 \rangle \leq 0
	\end{align*}
	which implies that $T' \lambda_2 = g$ holds $\mu_2$-almost surely on $[2b(0), \infty)$ and $T\mu_2 = 1$ holds $\lambda_2$-almost surely on $\Omega$. We summarise our findings.
	
	\begin{theo}\label{ThmSD2}
		 There exists a sequence $(\varepsilon_n)_{n \in \nn}$ and measures $\lambda_2 \in \M^+(\Omega)$ and $\mu_2 \in \M^+\big([2b(0), \infty)\big)$ such that $\lambda_{2, \varepsilon_n} \rightarrow \lambda_2$ and $\mu_{2, \varepsilon_n} \rightarrow \mu_2$ vaguely. The measures $\mu_2$ and $\lambda_2$ are \eqref{PROGRAM_P_2}- and \eqref{PROGRAM_D_2}-admissible, respectively, strong duality holds, i.e., $d_2 = p_2$ and $\mu_2$ is \eqref{PROGRAM_P_2}-optimal.
		
		In addition, if $\| \lambda_{2, \varepsilon_n} \| \rightarrow \| \lambda_2 \|$,
		then $\lambda_2$ is \eqref{PROGRAM_D_2}-optimal and the following complementary slackness conditions are satisfied: 
			\begin{align*}
				\int r_\theta\big(t,b(t)\big)\mu_2(d\theta)&= 1  &&\hspace{-1.5cm}\text{for $\lambda_2$-almost all\ } t\in (0,t_0],\\
				\int r_\theta\big(t,b(t)\big)\lambda_2(dt)&=r_\theta(t_0,x_0) &&\hspace{-1.5cm} \text{for $\mu_2$-almost all\ } \theta\in [2b(0),\infty).
			\end{align*}

	\end{theo}
	
	\begin{remark}
		Note that in this second set-up, we have strong duality and optimality of $\mu_2$ without any assumption, while in Theorem \ref{ThmSD1} we needed the additional assumption $\| \lambda_{1, \varepsilon_n} \| \rightarrow \| \lambda_1 \|$ in order to obtain strong duality and optimality of $\mu_1$. In both cases, the additional assumption is necessary for the optimality of $\lambda_1$ and $\lambda_2$, respectively, as well as for the complementary slackness conditions to hold. This is due to the fact that any \eqref{PROGRAM_D_2_eps}-feasible $\lambda$ is also \eqref{PROGRAM_D_2}-admissible as we consider $\lambda$ to be the null measure outside of $\Omega_\varepsilon$. However, a \eqref{PROGRAM_D_1_eps}-admissible $\mu$ does not have to be \eqref{PROGRAM_D_1_short}-admissible, as $\mu$ may not fulfil $T\mu = 1$ outside of $\Omega_\varepsilon$.
		\end{remark}

	\section{On the existence of representing measures}\label{SecExisRepMeasure}
	In this section, we investigate the existence of a representing measure. More precisely, given an analytic, concave boundary $b$ with $b(0)>0$ we want to prove the existence of a measure $\mu$ such that for all $t \in(0,t_0]$
	\begin{align*}
		1 = r\big(t, b(t)\big)=\int_{[2b(0),\infty)}r_\theta\big(t,b(t)\big)\mu(d\theta).
	\end{align*}

	 Now we see why the approach via linear programs is so fruitful for the question of representability: The complementary slackness conditions in Theorem \ref{ThmSD1} already state that
	\begin{align*}
		r\big(t, b(t)\big)=\int_{[2b(0),\infty)}r_\theta\big(t,b(t)\big)\mu_1(d\theta) = 1 \qquad  \text{ for $\lambda_1$-almost all\ } t\in (0,t_0],
	\end{align*}
	where $\mu_1$ and $\lambda_1$ are optimal measures for \eqref{PROGRAM_D_1_short} and \eqref{PROGRAM_P_1_short}, respectively. The analogous result holds for $\mu_2$ and $\lambda_2$ by Theorem \ref{ThmSD2}. Thus, if $\lambda_1$ or $\lambda_2$ put mass everywhere on $\Omega = (0, t_0]$, then $b$ is representable on $[0,t_0]$. 
	
	As we want to make use of the complementary slackness conditions as well as of strong duality throughout this section, we assume 
	\begin{assumption2}
	Let $b\colon[0,\infty)\to \rr$ satisfy Assumption (A) and let the assumptions of Theorems \ref{ThmSD1} and \ref{ThmSD2} hold true, in particular,  $\big\| \lambda_{1, \varepsilon^{(1)}_n} \big\| \to \| \lambda_1 \|$, $\big\| \lambda_{2, {\varepsilon^{(2)}_n}} \big\| \to \| \lambda_2 \|$ as $n\to\infty$, where $(\eps^{(i)}_n)_{n\in\nn}$ is a sequence   such that  $\lambda_{i,\eps^{(i)}_n}$  is an optimisers in $(P_{i,\eps_n})$ and $\lambda_i$ the vague limit of $\big(\lambda_{i,\eps^	{(i)}_n}\big)_{n\in\nn}$, $i=1,2$.
\end{assumption2}

		First we show that there exists a measure $\bar{\lambda}$ which is admissible in \eqref{PROGRAM_P_1_short} as well as in \eqref{PROGRAM_D_2} and already fulfils all constraints in both programs with equality.
	\begin{lemma}\label{lem:last_hitting}
		Let $\sigma_b$ be the last passage time of the Brownian motion $W$ to $b$ before time $t_0$, i.e.,
		\begin{align*}
			\sigma_b=\sup\{t\in[0,t_0]\colon W_t=b(t)\},
			\end{align*}
			where we set $\sigma_b=t_0$ on $\{W_t<b(t) \text{ for all } t\in[0,t_0]\}$, 
			 and define
		\begin{align*}
			\bar\lambda(dt)=\P(\sigma_b\in dt|W_0=0,\,W_{t_0}=x_0).
		\end{align*}
		Then, $\bar\lambda$ is admissible for \eqref{PROGRAM_P_1_short} as well as for \eqref{PROGRAM_D_2} and 
		\begin{align}
			\int r_\theta\big(t,b(t)\big)\bar\lambda(dt)= r_\theta(t_0,x_0) \qquad \text{ for any } \theta\in [2b(0),\infty).\label{eq:representable}
		\end{align}
	\end{lemma}
	\begin{proof}
		First, note that $\bar{\lambda} \in \M^+\big((0, t_0]\big)$. 
	 Let $\theta\geq 2b(0)$ and consider the additional linear boundary 
	 \begin{align*} 
		\ell_\theta(t)=\frac \theta 2+m_\theta\,t, 
		\end{align*}
		where $m_\theta\geq b^\prime(0)$. Note that $\ell_\theta$ lies above the boundary $b$. Let $\tau_\theta$ be the first passage time of $W$ to $\ell_\theta$. Then by Remark~\ref{rem:linear} the boundary $\ell_\theta$ is representable by $\mu_{\ell_\theta}=\exp\left(-\theta\:\!m_\theta\right)\delta_\theta$ and by Proposition \ref{prop:meth_images} it holds for   $t\in(0,t_0]$ and $ x< \ell_\theta(t)$ that 
		\begin{align*}
		\P(\tau_\theta\leq t\:\!|\:\! W_0=0,W_{t}=x)=\P(\tau_\theta<t \:\!|\:\!W_t=x)=\exp(-\theta\:\! m_\theta)\:\!r_\theta(t,x). 
		\end{align*}
		To show \eqref{eq:representable}, denote by $\P_{(0,0)}^{(s,x)}$ the law of a Brownian bridge from $(s, x)$ to $(0, 0)$. Here, we consider the Brownian bridge from $(s, x)$ to $(0, 0)$   as a process in reversed time. In particular, for a Brownian bridge from $(t_0,x_0)$ to $(0,0)$ note that the first passage time to $b$ coincides with the last passage time $\sigma_b$  of $b$ on $(0,t_0]$ for a Brownian bridge from $(0,0)$ to $(t_0,x_0)$. Using the strong Markov property of the Brownian bridge gives
		\begin{align*}
			r_\theta(t_0,x_0)&=\exp(\theta\:\! m_\theta)\,\P(\tau_\theta\leq t_0\:\!|\:\!W_0=0,W_{t_0}=x_0)\\
			&=\exp(\theta\:\!m_\theta)\,\!\int_{(0,t_0]} \P_{(0,0)}^{(t,b(t))}(\tau_\theta\leq t)\,\P_{(0,0)}^{(t_0,x_0)}(\sigma_b\in dt)\\
			&=\exp(\theta\:\!m_\theta)\,\int_{(0,t_0]} \exp(-\theta\:\!m_\theta)\, r_\theta(t,b(t))\,\P_{(0,0)}^{(t_0,x_0)}(\sigma_b\in dt)\\
			&=\int_{(0,t_0]} r_\theta(t,b(t)) \bar{\lambda}(dt).
		\end{align*}
		Hence, \eqref{eq:representable} holds and in particular, $\bar{\lambda}$ is admissible in both \eqref{PROGRAM_P_1_short} and \eqref{PROGRAM_D_2}.
	\end{proof}
	\begin{remark} 
		Note that $\bar{\lambda}$ is a strong contender to be optimal in both \eqref{PROGRAM_P_1_short} and \eqref{PROGRAM_D_2} as it already fulfils the constraint not only with inequality but even with equality everywhere. If $\bar{\lambda}$ really is optimal in both programs, then this yields a stochastic interpretation of \eqref{PROGRAM_P_1_short} and \eqref{PROGRAM_D_2} as the programs having the (conditional) last passage time distribution of a standard Brownian motion to $b$ as optimisers. This would be a nice symmetry with the programs \eqref{PROGRAM_D_1_short} and \eqref{PROGRAM_P_2}, where in both programs an optimiser is given by the (conditional) first passage time distribution. 
		\end{remark}
	
	Recall that $p_1,d_1, p_2$ and $d_2$ denote the optimal values of \eqref{PROGRAM_P_1_short}, \eqref{PROGRAM_D_1_short}, \eqref{PROGRAM_P_2} and \eqref{PROGRAM_D_2}, respectively. If we assume strong duality in both set-ups and since $\bar{\lambda}$ is admissible in both \eqref{PROGRAM_P_1_short} and \eqref{PROGRAM_D_2} by Lemma \ref{lem:last_hitting} we obtain
	\begin{align}\label{IneqLambdaBar}
		d_1 = p_1 \leq \big\| \bar{\lambda} \big\| \leq d_2 = p_2.
	\end{align}
	With the help of \eqref{IneqLambdaBar} we can provide sufficient conditions for a general concave, analytic $b$ to be representable.
	\newpage 
	
	\begin{theo}\label{ThmRepEqual}
		Assume that Assumption (B) is satisfied.  If either
		\begin{enumerate}[(i)]
			\item $d_1 = p_2$ or 
			\item $p_1 = d_2$, 
		\end{enumerate}
		then, $b$ is representable on $[0,t_0]$.
	\end{theo}
	\begin{proof}
		The proof follows directly from \eqref{IneqLambdaBar}:  If one of the conditions is fulfilled, we immediately have that $\bar{\lambda}$ is the optimal measure in both \eqref{PROGRAM_P_1_short} and \eqref{PROGRAM_D_2}. As $\bar{\lambda}$ puts mass everywhere in $(0, t_0]$, we conclude from either Theorem \ref{ThmSD1} or Theorem \ref{ThmSD2} that $r\big(t, b(t)\big) = 1$ for all $t \in (0, t_0]$ and so $b$ is representable on $[0,t_0]$.
	\end{proof}
	\begin{remark}
	 Of course, the conditions $d_1=p_2$ or $p_1=d_2$  in Theorem \ref{ThmRepEqual} are equivalent by \eqref{IneqLambdaBar}. We state the theorem in this way to stress that if the optimal values of the ``$\mu$-problems'' \eqref{PROGRAM_D_1_short} and \eqref{PROGRAM_P_2} or the optimal values of the ``$\lambda$-problems'' \eqref{PROGRAM_P_1_short} and \eqref{PROGRAM_D_2} agree, then we have representability.
		\end{remark}
	
	\begin{remark}
		The conditions $d_1=p_2$ or $p_1=d_2$  from Theorem \ref{ThmRepEqual} are a substantive improvement over the conditions that we usually impose to guarantee that $b$ is representable. Usually we would have to prove $r\big(t,b(t)\big) =1$ for all $t \in (0, t_0]$, now it suffices to show $d_1 = p_2$ or $p_1 = d_2$. In particular, these conditions can easily be checked in implementations (see Section \ref{SecNumStuRep} below).
			\end{remark}

	Moving away from the measure $\bar{\lambda}$, we can also give more sufficient conditions such that $b$ is representable on $[0,t_0]$. To this end, we first derive some properties of admissible measures in \eqref{PROGRAM_P_1_short} and \eqref{PROGRAM_D_2}, respectively.
	\begin{lemma}\label{lem:mass_t_0}
		Let $\lambda$ be admissible for \eqref{PROGRAM_P_1_short}. Then, for every $\varepsilon\in(0,t_0)$ we have  $\lambda([t_0-\varepsilon, t_0]) > 0$.
	\end{lemma}
	\begin{proof}
		Let $\lambda$ be admissible in \eqref{PROGRAM_P_1_short}. In particular,
		\begin{align*}
			\int_\Omega r_\theta\big(t,b(t)\big)\lambda(dt)\geq r_\theta(t_0,x_0) \qquad \text{ for any } \theta\in [2b(0),\infty),
		\end{align*}
		which implies
		\begin{align*}
			\int \exp\left(-\frac{1}{2}\theta^2\left(\frac{1}{t}-\frac{1}{t_0}\right)+\theta\left(\frac{b(t)}{t}-\frac{x_0}{t_0}\right) \right)\lambda(dt)\geq 1  \qquad \text{ for any } \theta\in [2b(0),\infty).
		\end{align*}
		Observe that we may write 
		\begin{align*}
			\int \exp\left(-\frac{1}{2}\theta^2\left(\frac{1}{t}-\frac{1}{t_0}\right)+\theta\left(\frac{b(t)}{t}-\frac{x_0}{t_0}\right) \right)\lambda(dt)= \int \alpha(t)\exp(-\beta(t)\big(\theta-\gamma(t)\big)^2)\lambda(dt)
		\end{align*}
		for certain positive functions $\alpha, \beta$ and $\gamma$ being bounded on $[0,t_0-\varepsilon]$ for each $\varepsilon\in(0,t_0)$. Now, if $\varepsilon\in(0,t_0)$ is such that $\lambda([t_0-\varepsilon,t_0])=0$, then dominated convergence yields the contradiction 
		\begin{align*}
			\int \exp\left(-\frac{1}{2}\theta^2\left(\frac{1}{t}-\frac{1}{t_0}\right)+\theta\left(\frac{b(t)}{t}-\frac{x_0}{t_0}\right) \right)\lambda(dt)\rightarrow0\mbox{ as } \theta \rightarrow \infty.
		\end{align*}
	\end{proof}
	\begin{lemma}\label{lem:no_mass_t_0}
		Let $\lambda$ be admissible for \eqref{PROGRAM_D_2}. Then, $\lambda(\{t_0\}) = 0$.
	\end{lemma}
	\begin{proof}
		Let $\lambda$ be admissible for \eqref{PROGRAM_D_2}. In particular, we have that for all $\theta \in [2b(0), \infty)$
		\begin{align*}
			\int_{(0, t_0]} \exp\left( -\frac{\theta^2}{2} \left( \frac{1}{t} -\frac{1}{t_0} \right) + \theta \left( \frac{b(t)}{t} - \frac{x_0}{t_0} \right) \right) \lambda (dt) \leq 1.
		\end{align*}
		Assume that $\lambda(\{t_0\}) > 0$, then we find
		\begin{align*}
			1 &\geq \int_{(0,t_0]}\exp\left( -\frac{\theta^2}{2} \left( \frac{1}{t} -\frac{1}{t_0} \right) + \theta \left( \frac{b(t)}{t} - \frac{x_0}{t_0} \right) \right) \lambda (dt)\\ &\geq
			 \exp\left( \theta \left( \frac{b(t_0)-x_0}{t_0} \right) \right) \cdot \lambda(\{t_0\}) \xrightarrow[\theta\to\infty]{}\ \infty
		\end{align*}
		as $b(t_0) > x_0$ by assumption. This is a contradiction. So $\lambda(\{t_0\}) = 0$.
	\end{proof}
	With this results we derive sufficient conditions guaranteeing that $b$ is representable.
	\begin{theo}\label{ThmRepSeq2}
		Let Assumption (B) be satisfied. Moreover, 
		assume that the integrability condition from Theorem \ref{ThmRepSeq} is satisfied, i.e., assume that there exist some $t^\ast>t_0$ such that
		\begin{align*}
			\int_{(0, \infty)} \exp\left(-\frac{\theta^2}{2 t^\ast} \right) \mu_i(d\theta) < \infty
		\end{align*}
		for $i =1, 2$, where  $\mu{_1}$ and $\mu_2$ are optimisers in \eqref{PROGRAM_D_1} and \eqref{PROGRAM_P_2} respectively.  Moreover, assume that one of the following conditions is met:
		\begin{enumerate}[(i)]
			\item $\lambda_1(\{t_0\}) = 0$,
			\item for every $\varepsilon\in(0,t_0)$ it holds $\lambda_2([t_0-\varepsilon, t_0]) > 0$, 
		\end{enumerate}
		where $\lambda_1$ and $\lambda_2$ are optimisers in \eqref{PROGRAM_P_1}  and \eqref{PROGRAM_D_2}, respectively. Then, $b$ is representable on $[0,t_0]$.
	\end{theo}
	\begin{proof} 
		First assume condition (i) holds. Lemma \ref{lem:mass_t_0} yields that $\lambda_1$ puts mass in every interval $[t_0-\varepsilon, t_0]$ for every $\varepsilon \in (0, t_0)$ but by assumption $\lambda_1(\{t_0\}) = 0$. Due to the complementary slackness conditions from Theorem \ref{ThmSD1}, we conclude that there exists a strictly increasing sequence $t_1, t_2, \ldots \nearrow t_0$ such that
		\begin{align*}
			\int r_\theta\big(t_n, b(t_n)\big) \mu_1(d\theta) = 1 \qquad \text{ for all } n \in \nn.
		\end{align*}
		Theorem \ref{ThmRepSeq} then implies that $b$ is representable on $[0,t_0]$. 
		
		If condition (ii) holds the proofs follows along the same lines with the help of Lemma \ref{lem:no_mass_t_0} and Theorem \ref{ThmSD2}.
	\end{proof}
	\begin{remark}
		Note that Theorem \ref{ThmRepSeq2} also offers sufficient conditions for $b$ to be representable that are easier to check than the usual condition $r\big(t, b(t)\big)=1$ for all $t\in (0,t_0]$. In Chapter \ref{ChapCompMeth}, we investigate a new method to obtain numerical candidates for $\mu_1$ and $\lambda_1$ or $\mu_2$ and $\lambda_2$, respectively. Then, it is rather straightforward to check whether these measures fulfil the conditions from Theorem~ \ref{ThmRepSeq2} at least numerically.
		\end{remark}
	
	\section{Computational method for the linear programming approach}\label{ChapCompMeth}
	After a short introduction to existing numerical approaches for the inverse method of images, we give a convergence result for discretised versions of our programs and based on that new algorithm. We also provide error bounds for the numerical distribution function of the first passage time to a boundary $b$. Moreover, we numerically investigate representability, i.e., we present a numerical study of the assumptions in Theorem \ref{ThmRepEqual}. 
	
	\subsection{Existing computational approaches for the inverse method of images}
	\sectionmark{Existing computational approaches}
	We will now discuss some approximation methods for the inverse method of images that have already been investigated. 
	 
	In \cite{LRD02}, 
 the authors consider 	for a given a boundary $b$, the equation $r_\mu\big(t,b(t)\big)=1$
	for $t \in (0,\infty)$, where $\mu$ has to be found in the set of signed measure. To approximate the actual but unknown 
	$\mu$ with a signed measure $\tilde{\mu}$ such that the boundary $\tilde{b}$ generated by $\tilde{\mu}$ is close to $b$, it is assumed that
	$		\mu=\sum_{r=1}^N w_r \:\!\delta_{\theta_r} $
		with unknown weights $w_r\in \rr$, $r = 1, \ldots, N$ for some $N \in \nn$. Moreover, they choose a set of increasing time points $t_s$, $s = 1, \ldots, 2N$, to obtain the simplified equations
	\begin{align}\label{eq:disctime}
		1 = r_\mu\big(t_s,b(t_s)\big) =\sum_{r=1}^N w_r \exp\left( \theta_r \frac{b(t_s)}{t_s} - \frac{\theta_r^2}{2t_s} \right), \quad s=1, \ldots, 2N.
	\end{align}
	If the values of $\theta_r$ are pre-assigned, only $N$ time points are required to solve \eqref{eq:disctime} which is linear in $w_r$. The idea behind the discretisation in \eqref{eq:disctime} is that the given boundary $b$ and the boundary $\tilde{b}$ generated by the signed measure $\tilde{\mu}$ through the method of images are equal at the time points $t_s$, $s=1,\ldots, 2N$.  
	The authors emphasize that a higher number of time points may be desirable to increase the accuracy of the approximation of $b$ with $\tilde{b}$ but this may of course turn the system \eqref{eq:disctime} singular.
	
	The authors investigate the approximation error by showing that for the first passage times $\tau_b$ and $\tau_{\tilde{b}}$ to the boundaries $b$ and $\tilde{b}$, respectively, it holds for all $t >0$ that
	\begin{align*}
		\big| \P(\tau_b < t) - \P(\tau_{\tilde{b}} < t) \big| \to 0 \quad \text{as }\ \bar{\varepsilon}_t \coloneqq \sup_{0<s<t} \big\vert \tilde{b}(s) -b(s) \big\vert \to 0.
	\end{align*}
		In other words: the distribution functions of $\tau_b$ and $\tau_{\tilde{b}}$ are close if $b$ and $\tilde{b}$ are close. However, no convergence result for the algorithm itself is provided. The investigated examples of square-root and parabolic boundaries show fairly good approximations with deviations of $\tilde{b}$ from $b$ especially for small $t$. 
	
	In \cite{ZIP13}, Zipkin refines the numerical methods developed by among others \cite{LRD02}.
		Given a boundary~$b$, \cite{ZIP13} uses the same discretisations as \cite{LRD02} but requires non-negative weights $w_r$. Here the number of points with non-negative mass and the number of time discretisations are denoted by $J$ and $I$, respectively, and do not satisfy $I=2J$ in general. Then, for $w = (w_j)_{j=1, \ldots, J}$ and  $\one=(1,\ldots,1)\in \R^I$ Equation \eqref{eq:disctime} can be restated as $Mw = \one$ for a suitable matrix $M$. \cite{ZIP13} relaxes this setting slightly by choosing slack variables $s\in \rr^I$ and by fixing a vector $p\in (0,\infty)^I$. Then, the following linear program has to be solved
	\begin{alignat*}{3}
		&\text{minimise} &\quad &p^Ts\\
		&\text{subject to} & &Mw +s = \one,\\
		& & & w \geq 0,\\
		& & & s \geq 0.
	\end{alignat*}
	The linear program's objective function minimises the $p$-weighted deviations $s$. Solving the linear program then give the weights $w_j, j=1, \ldots, J$, such that $\tilde{b}$ generated by $\tilde{\mu} = \sum_j w_j \delta_{\theta_j}$ most closely resembles $b$. \cite{ZIP13} does not give any convergence results for his algorithm. However, he finds in examples that the approximation $\tilde{b}$ of $b$ generally works well on the discretised time interval $[t_1, t_I]$ but deviates from $b$ outside that interval. The algorithm can nevertheless be seen as a substantial improvement on the algorithm from \cite{LRD02}. \cite{ZIP13} also presents a generalisation of how to include not only point measures but also more general positive, $\sigma$-finite measures with densities.
	
	Although the algorithms set out by \cite{LRD02} and \cite{ZIP13} provide good approximations, they suffer from some drawbacks: in both cases, the measure has to be discretised as a weighted sum of point measures as well as the time axis has to be discretised. Moreover, for neither algorithm a convergence result is provided. Both problems are addressed in the following section.
	
	\subsection{Convergence results and a new algorithm}\label{SecConvRes}
	
	In this section, we present an algorithm for the inverse method of images and provide convergence results. This section is inspired by methods for American options set out in \cite{LEN17,CHR14}. 
	
	Let us start by discretising the linear problem. To this end, let $\mu_i \in \M^+\big([2b(0), \infty)\big)$, $i \in \nn$,  and for $n\in\nn$ let 
	\begin{align*}
		U_n \defined = \left\{ \sum_{i=1}^n a_i\mu_i \,\colon a \in [0,\infty)^n\right\}
	\end{align*}
	 denote the positive cone generated by the measures $\mu_i, i=1, \ldots, n$, and denote by $U_\infty$ the closure of $\bigcup_{n \in \nn} U_n$ with respect to the vague topology.  Restricting the linear program \eqref{PROGRAM_P_2} to measures in  $U_n$ results in
	\begin{equation*}
		\label{PROGRAM_P_2_n}
		\boxed{\begin{aligned}
			& \text{minimise} & & \int r_\theta(t_0,x_0)\mu(d\theta) \\
			& \text{subject to} & & \mu \in U_n,\\
			&&& \int r_\theta\big(t,b(t)\big)\mu(d\theta)\geq 1 \qquad \text{for any } t \in (0,t_0].
		\end{aligned}\tag{$P_{2,n}$}
	}
	\end{equation*}

	Now we prove the following existence and consistency result for our simplified program~\eqref{PROGRAM_P_2_n}.
	\begin{propo}\label{PropConv1}
			Assume that there exists $C \in (0,\infty)$ such that $C \cdot \mu_1$ is admissible in $(P_{2,1})$ and therefore in any \eqref{PROGRAM_P_2_n} for $n \in \nn \cup \{\infty\}$. 
		\begin{enumerate}[(a)]
			\item For $n \in \nn \cup \{ \infty \}$ the optimal value $p_{2, n}$ of the linear program \eqref{PROGRAM_P_2_n} is attained at some admissible measure $\mu_{2, n}^\ast$ and satisfies $p_{2, n} \leq C \Vert \mu_1 \Vert$. Moreover, for $m \leq n$ the measure $\mu_{2,m}^\ast$ is \eqref{PROGRAM_P_2_n}-admissible and it holds that $p_{2, m} \geq p_{2, n}$.
			\item\label{PropConv1_item_2} There exists a subsequence of optimisers $(\mu_{2, n_k}^\ast)_{k\in\nn}$ and a $(P_{2,\infty})$-admissible measure $\nu_\infty$ such that $\mu_{2, n_k}^\ast \to \nu_\infty$ vaguely as $k\to\infty$. Moreover,
			\begin{align*}
				p_{2, \infty} \leq \int r_\theta(t_0, x_0) \nu_\infty (d\theta) \leq \inf_{n \in \nn} p_{2, n} = \lim_{n \rightarrow \infty} p_{2, n}.
			\end{align*}
			\item If there exists a sequence $(\xi_n)_{n\in\nn}$ with $\xi_n \in U_n$ and   $\sup_{n\in\nn} \|\xi_n\|<\infty$ such that $(\xi_n)_{n\in\nn}$ converges vaguely to some $(P_{2,\infty})$-optimal measure $\mu_\infty^\ast$ as $n \rightarrow \infty$ and if
			\begin{align*}
				\lim_{n \rightarrow \infty}\, \sup_{t \in (0, t_0]} \frac{\left\vert \int r_\theta\big(t, b(t)\big) \mu_\infty^\ast (d\theta) - \int r_\theta\big(t, b(t)\big) \xi_n (d\theta) \right\vert}{\int r_\theta \big(t, b(t)\big) \mu_1(d\theta)} = 0,
			\end{align*}
			then $\nu_\infty$ from (b) is $(P_{2,\infty})$-optimal. 
			Moreover, we have
			\begin{align*}
				p_{2, \infty} = \int r_\theta(t_0, x_0) \nu_\infty (d\theta) = \inf_{n \in \nn} p_{2, n} = \lim_{n \rightarrow \infty} p_{2, n}.
			\end{align*}
		\end{enumerate}
	\end{propo}
	\begin{proof}
		\begin{enumerate}[(a)]
			\item 			
				For any $n \in \nn\cup \{\infty\}$ we reformulate the program \eqref{PROGRAM_P_2_n} by absorbing $r_\theta(t_0, x_0)$ into~$\mu$. The linear program now reads
			\begin{align}
				\label{PROGRAM_P_2_n'}
				\boxed{
				\begin{aligned}
					& \text{\normalfont{minimise}} & & \Vert \mu \Vert \\
					& \text{\normalfont{subject to}} & & \mu \in U_n,\\
					&&& \int \frac{r_\theta\big(t,b(t)\big)}{r_\theta(t_0, x_0)} \mu(d\theta) \geq 1 \qquad \text{for any } t \in (0,t_0].
				\end{aligned} \tag{$P_{2,n}'$}
			}
			\end{align}
			Note that $\mu$ is admissible in \eqref{PROGRAM_P_2_n} if and only if ${\mu'}$ given by 
			\begin{align*}
				\frac{d\mu'}{d\mu}(\cdot)=r_\cdot(t_0,x_0)
			\end{align*}
			is admissible in \eqref{PROGRAM_P_2_n'}. 
			
			Let $\A_{P'_{2,n}}$ be the set of all admissible measures in \eqref{PROGRAM_P_2_n'}. 
			Since $\mu_a \defined = C \cdot \mu_1 \in U_n$ is \eqref{PROGRAM_P_2_n}-admissible for every $n \in \nn \cup \{\infty\}$ by assumption, the measure $\mu'_a$ defined by
			\begin{align*}
			\frac{d	\mu'_a}{d \mu_a} (\cdot)\defined = r_\cdot (t_0, x_0) 
			\end{align*}
			 is \eqref{PROGRAM_P_2_n'}-admissible. Moreover, every potential minimiser $\mu^\ast \in \A_{P'_{2,n}}$ satisfies
			\begin{align*}
				\Vert \mu^\ast \Vert \leq \Vert \mu'_a \Vert \eqqcolon \rho,
			\end{align*}
			where $\rho < \infty$. In particular, any potential solution $\mu^\ast \in \A_{P'_{2,n}}$ is contained in the vaguely compact ball $B_\M(\rho) = \left\{ \mu \in \M\big([2b(0), \infty)\big) \colon \Vert \mu \Vert \leq \rho \right\}$. Therefore, in the linear program \eqref{PROGRAM_P_2_n'} it is sufficient to consider measure from the set
			\begin{align*}
				\cA_{P'_{2,n}} \cap B_\M(\rho) = \bigcap_{t \in (0, t_0]} H(t) \cap U_n \cap B_\M(\rho),
			\end{align*}
			where 
			\begin{align*} 
				H(t) \defined = \left\{ \mu \in \M(\rr) \colon \int \frac{r_\theta\big(t,b(t)\big)}{r_\theta(t_0, x_0)} \mu(d\theta) \geq 1 \right\}.
				\end{align*}
				\newpage 
				 Since $H(t)$ and $U_n$ are closed with respect to the vague topology,  and $B_\M$  is vaguely compact due to the Alaoglu-Bourbaki theorem, cf.\ \cite[Theorem 23.5]{MV92}, also $\A_{P'_{2,n}} \cap B_\M(\rho)$ is vaguely compact. Then, the optimal value $p'_{2,n}$ of the linear program \eqref{PROGRAM_P_2_n'} is attained by some measure $\mu'_{2,n} \in \A_{P'_{2,n}} \cap B_\M(\rho)$ as $\mu \mapsto \Vert \mu \Vert$ is lower semi-continuous with respect to the vague topology. Then, the optimal value $p_{2,n}$ in \eqref{PROGRAM_P_2_n} is attained by $\mu^\ast_{2,n}$ where 
				 \begin{align*}
				 	\frac{d\mu^\ast_{2,n}}{  d\mu'_{2,n} }(\cdot)= \frac{1}{r_\cdot(t_0, x_0)} 
				 	\end{align*}
				 	 and it holds that 
				 \begin{align*}
				 	p_{2,n}=\int r_\theta(t_0,x_0)\:\!\mu^\ast_{2,n}(d\theta)=\|\mu_{2,n}^\prime\|\leq \rho=C\|\mu_1\|.
				 \end{align*}
			\item Since $\Vert \mu'_{2,n} \Vert \leq \rho$ for any $n \in \nn$ and as the total variation unit ball is vaguely compact, there exists a subsequence $(\mu_{2,n_k}')_{k\in\nn}$ and a measure $ \nu'_\infty\in U_\infty\cap B_\M(\rho)$  such that $\mu_{2,n_k}' \to \nu_\infty'$ vaguely as $k\to\infty$. Then, we obtain that the subsequence $\big(\mu^\ast_{2,n_k} \big)_{k\in\nn}$ 
			converges vaguely to $\nu_\infty$, where 
			\begin{align*}
				\frac{d \nu_\infty}{d\nu_\infty'}(\cdot)=\frac{1}{r_\cdot(t_0,x_0)}.
				\end{align*}
			Let  $t \in (0, t_0)$. Then, the mapping $[2b(0), \infty) \to \rr, \theta \mapsto \frac{r_\theta(t,b(t))}{r_\theta(t_0, x_0)}$ vanishes at infinity (which does not hold true for $t=t_0$), so \cite[Theorem 30.6]{BAU01} yields
			\begin{align*}
				\int \frac{r_\theta\big(t,b(t)\big)}{r_\theta(t_0, x_0)} \:\!\nu_\infty'(d\theta) = \lim_{k\rightarrow \infty} \int \frac{r_\theta\big(t,b(t)\big)}{r_\theta(t_0, x_0)}\:\! \mu_{2,n_k}'(d\theta) \geq 1
			\end{align*}
			as $\mu'_{2,n_k}$ is ($P'_{2,n_k}$)-admissible and hence, 
			\begin{align*}
				T\nu_\infty(t)=\int r_\theta\big(t,b(t)\big)\nu_\infty(d\theta)=\int \frac{r_\theta\big(t,b(t)\big)}{r_\theta(t_0, x_0)} \:\!\nu_\infty'(d\theta)\geq 1 \qquad \text{for all $t\in(0,t_0)$.}
			\end{align*}
			Since $T\nu_\infty$ is continuous, we conclude that
			\begin{align*} 
				T\nu_\infty(t_0)=\int r_\theta\big(t_0,b(t_0)\big)\nu_\infty(d\theta)\geq 1.
			\end{align*}
			
		Thus, $\nu_\infty$ is $(P_{2,\infty})$-admissible. In particular,  $p_{2,\infty} \leq \int r_\theta(t_0, x_0) \nu_\infty (d\theta)$. Moreover, using vague convergence and \cite[Lemma 30.3]{BAU01} we find that $\Vert \nu_\infty' \Vert \leq \liminf_{k \rightarrow \infty} \Vert \mu_{2, n_k}' \Vert$. As $(p_{2,n})_{n\in\nn}$ is monotone, we conclude that
			\begin{align*}
				p_{2, \infty} \leq \int r_\theta(t_0, x_0) \:\!\nu_\infty (d\theta) = \Vert \nu_\infty' \Vert \leq \liminf_{k \rightarrow \infty} \Vert \mu_{2, n_k}' \Vert = \liminf_{k\to\infty} p_{2,n_k} = \lim_{n \rightarrow \infty} p_{2, n}.
			\end{align*}
			\item Let $n \in \nn$ and define $\eta_n \defined = \xi_n + \varepsilon_n \mu_1$, where
			\begin{align*}
				\varepsilon_n \defined = \sup_{t \in (0, t_0]} \frac{\left\vert \int r_\theta\big(t, b(t)\big) \mu_\infty^\ast (d\theta) - \int r_\theta\big(t, b(t)\big) \xi_n (d\theta) \right\vert}{\int r_\theta \big(t, b(t)\big) \mu_1(d\theta)}.
			\end{align*}
			Since $\mu_{\infty}^\ast$ is $(P_{2,\infty})$-optimal by assumption and therefore $(P_{2,\infty})$-admissible, we find for all $t \in (0, t_0]$ that
			\begin{align*}
				\int r_\theta \big(t, b(t)\big) \:\!\eta_n(d\theta) -1
				\geq \, &\int r_\theta \big(t, b(t)\big) \:\!\eta_n(d\theta) - \int r_\theta \big(t, b(t)\big) \mu_\infty^\ast(d\theta)\\
				= \, &\int r_\theta \big(t, b(t)\big)\:\! \mu_1(d\theta) \left(\!\varepsilon_n - \frac{\int r_\theta \big(t, b(t)\big)\:\! \mu_\infty^\ast(d\theta) - \int r_\theta \big(t, b(t)\big)\:\! \xi_n(d\theta)}{\int r_\theta \big(t, b(t)\big) \:\!\mu_1(d\theta)} \!\right)\\
				\geq \, &0, 
			\end{align*}
			which means that $\eta_n$ is \eqref{PROGRAM_P_2_n}-admissible. Using \ref{PropConv1_item_2} we obtain
			\begin{align}
				\begin{split}
				p_{2,\infty} \leq \int r_\theta (t_0, x_0) \nu_\infty(d\theta) \leq p_{2,n} &\leq \int r_\theta (t_0, x_0) \eta_n(d\theta)\\
				&= \int r_\theta (t_0, x_0) \xi_n(d\theta) + \varepsilon_n \int r_\theta (t_0, x_0) \mu_1(d\theta).
				\end{split}\label{eq:PropConv1_1}
						\end{align}
			Since  $(\xi_n)_{n\in\nn}$ converges  vaguely to $\mu_{\infty}^\ast$ and  $\sup_{n\in\nn}\|\xi_n\|<\infty$ by assumption, we conclude from \cite[Theorem 30.6]{BAU01} that 
			\begin{align}
				\lim_{n \rightarrow \infty} \int r_\theta (t_0, x_0) \xi_n(d\theta) = \int r_\theta (t_0, x_0) \mu_{\infty}^\ast (d\theta) = p_{2,\infty}.\label{eq:PropConv1_2}
			\end{align}
			Combining \eqref{eq:PropConv1_1}, \eqref{eq:PropConv1_2} and $\lim_{n \rightarrow \infty} \varepsilon_n = 0$  yields
			\begin{align*}
				p_{2,\infty} \leq \int r_\theta (t_0, x_0) \nu_\infty(d\theta) \leq \lim_{n \rightarrow \infty} p_{2,n} \leq \lim_{n \rightarrow \infty} \int r_\theta (t_0, x_0) \eta_n(d\theta) = p_{2,\infty}.
			\end{align*}
			In particular, $\nu_\infty$ is $(P_{2,\infty})$-optimal. 
 \end{enumerate}
	\end{proof}
		An analogous statement holds for \eqref{PROGRAM_D_1}.
\begin{propo}\label{PropConv2}
		 For $n \in \nn \cup \{\infty\}$ consider the linear program
		\begin{align}\label{PROGRAM_D_1_n}
			\boxed{\begin{aligned}
				& \text{\normalfont{maximise}} & & \int r_\theta(t_0,x_0)\mu(d\theta) \\
				& \text{\normalfont{subject to}} & & \mu \in U_n,\\
				&&& \int r_\theta\big(t,b(t)\big)\mu(d\theta)\leq 1 \qquad  \text{for any } t \in (0,t_0]
			\end{aligned} \tag{$D_{1,n}$}
		}
		\end{align}
		and assume that there exists $C \in (0,\infty)$ such that $C \cdot \mu_1$ is admissible in $(D_{1,1})$.
		\begin{enumerate}[(a)]
			\item For $n \in \nn \cup \{ \infty \}$, the optimal value $d_{1, n}$ of the linear program \eqref{PROGRAM_D_1_n}  is attained at some admissible measure $\mu_{1, n}^\ast$. Moreover, it holds that $d_{1, m} \leq d_{1, n}\leq 1$ for $m \leq n$.
			\item There exists a subsequence of optimisers $(\mu_{1, n_k}^\ast)_{k\in\nn}$ and a $(D_{1,\infty})$-admissible measure $\nu_{\infty}$ such that $\mu_{1, n_k}^\ast \to \nu_{\infty}$ vaguely as $k\to\infty$. Moreover,
			\begin{align*}
				d_{1, \infty} \geq \int r_\theta(t_0, x_0)\:\! \nu_{\infty} (d\theta) \geq \sup_{n \in \nn} d_{1, n} = \lim_{n \rightarrow \infty} d_{1, n}.
			\end{align*}
			\item If there exists a sequence $(\xi_n)_{n\in\nn}$  with $\xi_n \in U_n$ and   $\sup_{n\in\nn} \|\xi_n\|<\infty$ such that $(\xi_n)_{n\in\nn}$ converges vaguely to some $(D_{1,\infty})$-optimal measure $\mu_\infty^\ast$ as $n \rightarrow \infty$ and if
			\begin{align*}
				\lim_{n \rightarrow \infty} \, \inf_{t \in (0, t_0]} \frac{\left\vert \int r_\theta\big(t, b(t)\big) \mu_\infty^\ast (d\theta) - \int r_\theta\big(t, b(t)\big) \xi_n (d\theta) \right\vert}{\int r_\theta \big(t, b(t)\big) \mu_1(d\theta)} = 0,
			\end{align*}
			then $\nu_\infty$ from (b) is $(D_{1,\infty})$-optimal.  
			Moreover, it holds that
			\begin{align*}
				d_{1, \infty} = \int r_\theta(t_0, x_0)\:\! \nu_\infty (d\theta) = \sup_{n \in \nn} d_{1, n} = \lim_{n \rightarrow \infty} d_{1, n}.
			\end{align*}
		\end{enumerate}
	\end{propo}
	\begin{proof}
		As in the proof of Lemma \ref{Lemma_exis_1_eps} consider the modified program
		\begin{align}
			\label{PROGRAM_D_1_n'}
			\boxed{\begin{aligned}
				& \text{\normalfont{maximise}} & & \langle \psi ,\mu \rangle \\
				& \text{\normalfont{subject to}} & & \mu \in U_n,\\
				&&& \int \frac{r_\theta\big(t,b(t)\big)}{r_\theta(t_0, b(t_0))} \mu(d\theta) \leq 1 \qquad \text{for any } t \in (0,t_0], 
			\end{aligned} \tag{$D_{1,n}'$}
		}
		\end{align}
		where $\psi\colon[2b(0),\infty)\to \R,$ $\theta\mapsto\exp\left(-\theta \frac{b(t_0)-x_0}{t_0}\right)$. By \eqref{eq:D_1_eps_mod_admissible_norm_bounded_by_1} any \eqref{PROGRAM_D_1_n'}-admissible $\mu$ satisfies $ \Vert \mu \Vert\leq 1$. 
		Thus, all admissible measures in \eqref{PROGRAM_D_1_n'} are contained in $B_\M(1) =\big \{ \mu \in \M\big([2b(0), \infty)\big) \colon \Vert \mu \Vert \leq 1 \big\}$. The rest of the proof follows quite similar to the proof of Proposition \ref{PropConv1}.
	\end{proof}
	Similar results can be formulated and proven for the programs \eqref{PROGRAM_P_1} and \eqref{PROGRAM_D_2}, for more details see Appendix~A.2 in \cite{DissOskar}.

	Observe that for $\mu_i=\delta_{\theta_i} $ with $\{\theta_i\}_{i\in\nn}$ a dense subset of $[2b(0), \infty)$ the set $U_\infty$   is vaguely dense in $\M^+\big([2b(0), \infty)\big)$ by \cite[Theorem 30.4]{BAU01}. Thus, we can restrict ourselves to the point measures  $\mu_i=\delta_{\theta_i} $ in Propositions \ref{PropConv1} and \ref{PropConv2}  and still are able to approximate every possible representing measure $\mu$ arbitrarily close.  The inclusion of measures with densities, as done in \cite{ZIP13}, is therefore not necessary for the convergence of our algorithm but may be useful for numerical reasons.
	
	For our algorithm we further simplify the discretised linear problem \eqref{PROGRAM_P_2_n} and only consider the interval $I=[2b(0), 2b(0) +l]$ for some  $l>0$. 
	 Let $\mu_i=\delta_{\theta_i}$ be the point measures in $\theta_i=2b(0)+\frac{i-1}{n}l$, $i=1, \ldots, n$, and let again $U_n \defined = \{ \sum_{i=1}^n a_i\mu_i \colon a \in [0,\infty)^n\}$  be the positive cone generated by the measures $\mu_i, i=1, \ldots, n$. The restriction \eqref{PROGRAM_P_2_n'} of the linear program \eqref{PROGRAM_P_2_n} to measures in $U_n$	 can then be simplified to
	\begin{align*}
	\boxed{
			\begin{aligned}
			& \text{minimise} & & \sum_{i=1}^n r_{\theta_i}(t_0,x_0) a_i \\
			& \text{subject to} & & a \in [0,\infty)^n,\\
			&&& \sum_{i=1}^n r_{\theta_i}\big(t,b(t)\big) a_i \geq 1 \qquad \text{for any } t \in (0,t_0].
		\end{aligned}
	}
	\end{align*}
	For the implementation we make use of the cutting plane algorithm described in \cite{LW92} and \cite{IW10}, where  also good convergence results are provided. We  now present our algorithm:
	
	\bigskip 
	
	\noindent $\boxed{\text{\textbf{Step 1:}}}$ Let the set of initial constraints $\Gamma_1 \subset (0,t_0]$ be the set $\{t_0\}$. Choose a maximum number of iterations $k_\mathrm{max}$ and set $k = 1$.
	
	\medskip		

	 \noindent $\boxed{\text{\textbf{Step 2:}}}$ Calculate a solution $a^{(k)} \in [0,\infty)^n$ of the finite dimensional linear program
		\begin{align*}
			\boxed{
			\begin{aligned}
				& \text{minimise} & & \sum_{i=1}^n r_{\theta_i}(t_0,x_0) a_i^{(k)} 
				 \\
				& \text{subject to} & & a^{(k)} \in [0,\infty)^n,\\
				&&& \sum_{i=1}^n r_{\theta_i}\big(t,b(t)\big) a_i^{(k)} 
				 \geq 1 \qquad \text{for all } t \in \Gamma_k.
			\end{aligned}
		}
		\end{align*}
		
			\medskip
		
		\noindent	$\boxed{\text{\textbf{Step 3:}}}$ Determine  a point $t^{(k)} \in (0, t_0]$, where the constraint is most severely violated, i.e.,
		\begin{align*}
			t^{(k)} \in \argmin_{t \in (0, t_0]} \left\{ \sum_{i=1}^n r_{\theta_i}\big(t,b(t)\big) a_i^{(k)}
			\right\}.
		\end{align*}
		
			\medskip
			
		\noindent $\boxed{\text{\textbf{Step 4:}}}$ Add the point $t^{(k)}$		to the set of constraints, i.e., set $\Gamma_{k+1} \defined = \Gamma_k \cup \big\{t^{(k)}\big\}$.
		
		\medskip
		
		\noindent $\boxed{\text{\textbf{Step 5:}}}$ If the maximum number of iterations is reached, i.e., if $k = k_\mathrm{max}$, terminate the algorithm and output the approximate solution
		\begin{align*}
			\tilde{a} \defined = a^{(k_\mathrm{max})}.
		\end{align*}
		Otherwise, increase the iteration counter $k \leadsto k+1$ and return to the second step.
		
		\bigskip 
	
\begin{remark}	\hfill
	\begin{enumerate}[(a)]
		\item 
	In Proposition \ref{PropConv1} we do not assume the measures $\mu_i, i \in \nn$, to be point measures as we do in our algorithm. This indicates that the algorithm also works for a much larger class of ``auxiliary measures'' $\mu_i$. 
	\item Our algorithm requires fewer assumptions than the algorithm proposed in \cite{ZIP13} as we only choose the points $(\theta_i)_{i=1,\ldots,n}$ in $[2b(0), \infty)$ where the algorithm can put point masses and the initial constraint at $t_0$. Then, our algorithm ``chooses'' the next points where the constraint is evaluated. In contrast, in \cite{ZIP13}, both the division of $[2b(0), \infty)$ as well as the points $t_1, \ldots, t_m$ where the constraint is evaluated have to be chosen in advance. 
\end{enumerate}
\end{remark}
	
	For the linear program \eqref{PROGRAM_D_1}, we use similar simplifications and an analogous algorithm.
	
	 For the linear programs \eqref{PROGRAM_D_2} and \eqref{PROGRAM_P_1}, i.e., for the ``$\lambda$''-problems, we divide the interval $(0, t_0]$ into $n_\lambda$ equidistant points $0<t_1 < \ldots < t_{n_\lambda} = t_0$ and let $\lambda_i=\delta_{t_i}$ for $i = 1, \ldots, n_\lambda$. We then use analogous algorithms which are initialised with $\Gamma_1 = \{2b(0)\}$.
	
	If we have obtained a candidate representing measure from one of the above outlined algorithms, it is interesting to know the difference of the cumulative distribution function generated by this measure and the actual cumulative distribution function. If a measure $\mu$ represents $b$ 
	on $(0, t_0]$, then 
	the cumulative distribution function $F$ of $\tau$ is given by
	\begin{align*}
		F(t)
		&= 1-\Phi\left( \frac{b(t)}{\sqrt{t}} \right) + \int_{(-\infty,b(t))} r_\mu(t,x)\, p_t(0,x) dx,
	\end{align*}
	recall \eqref{eq:density}.
	If we now have obtained a measure $\tilde{\mu}$ from our algorithms, we define $\tilde{r}(t,x) = r_{\tilde\mu}(t,x)= \int_{[2b(0), \infty)} r_\theta(t,x) \tilde{\mu}(d \theta)$ and approximate the true cumulative distribution function $F$ by
	\begin{align*}
		\widetilde{F}(t) = 1-\Phi\left( \frac{b(t)}{\sqrt{t}} \right) + \int_{(-\infty,b(t))} \tilde{r}(t,x)\, p_t(0,x) dx.
	\end{align*}
	\begin{propo}\label{PropNumBound}
		Let $F$ be the true first passage time distribution to a boundary $b$ and $\tilde{\mu}$ a measure such that for all $t \in (0, t_0]$
		\begin{align*}
			1-\zeta_1 \leq \frac{1}{\tilde{r}_{\tilde\mu}\big(t, b(t)\big)} \leq 1+\zeta_2
		\end{align*}
		for some $\zeta_1 \in [0,1)$ and $\zeta_2 >0$. Then,
		\begin{align*}
			\sup_{t \in (0, t_0]} \left\vert F(t) - \widetilde{F}(t) \right\vert \leq \max \left( \zeta_1, \zeta_2 \right).
		\end{align*}
		In particular, we have
		\begin{align*}
			\lim_{\substack{\zeta_1 \searrow 0\\ \zeta_2 \searrow 0}} \, \, \sup_{t \in (0, t_0]}\Big \vert F(t) - \widetilde{F}(t)\Big \vert = 0.
		\end{align*}
	\end{propo}
	\begin{proof}
		Since by assumption, we have for all $t \in (0, t_0]$ that
		\begin{align*}
			1-\zeta_1 \leq \frac{1}{\tilde{r}_{\tilde\mu}\big(t, b(t)\big)} \leq 1+\zeta_2
		\end{align*}
		for some $\zeta_1 \in [0,1)$ and $\zeta_2 >0$, Proposition \ref{prop:meth_images} yields for $t \in (0, t_0]$ and $x < b(t)$
		\begin{align*}
			(1-\zeta_1)\:\!\tilde{r}_{\tilde\mu}(t,x)\leq \P(\tau\leq t\:\!|\:\!W_{t}=x)\leq (1+\zeta_2)\:\!\tilde{r}_{\tilde\mu}(t,x).
		\end{align*}
		 Hence, using \eqref{eq:density} we obtain
		\begin{align*}
			F(t) &= 1-\Phi\left( \frac{b(t)}{\sqrt{t}} \right) + \int_{(-\infty,b(t))} \P(\tau \leq t\:\! \vert\:\! W_t = x) \:\! p_t(0,x) dx\\
			&\leq 1-\Phi\left( \frac{b(t)}{\sqrt{t}} \right) + \int_{(-\infty,b(t))} (1+\zeta_2) \:\!\tilde{r}(t,x)\:\! p_t(0,x) dx
					\leq (1+\zeta_2) \:\!\widetilde{F}(t).
		\end{align*}
		 In particular, we  conclude
		\begin{align*}
			F(t) - \widetilde{F}(t) \leq \zeta_2.
		\end{align*}
		Analogously, we find
		\begin{align*}
			F(t) - \widetilde{F}(t) \geq -\zeta_1.
		\end{align*}
		To summarise, it holds that
		\begin{align*}
			\sup_{t \in (0, t_0]} \vert F(t) - \widetilde{F}(t) \vert \leq \max \left( \zeta_1, \zeta_2 \right).
		\end{align*}
	\end{proof}

	\subsection{Numerical study of representability}\label{SecNumStuRep}
	
	In this section, we present numerical results for our algorithm for three concave, analytic boundaries $b$  and one convex, analytic boundary $b$ and investigate their representability. 
	For all boundaries, we choose 
	\begin{itemize}
		\item $t_0 = 1$,
		\item $x_0 = b(t_0)-1$ and 
		\item the maximum number of iterations $k_\mathrm{max} = 20$.
		\end{itemize}
		 Moreover, for the ``$\mu$-problems'' \eqref{PROGRAM_D_1_short} and \eqref{PROGRAM_P_2} we select a length of $l=5$ and discretise the interval $[2b(0), \infty)$ by setting $n=100$ equidistant points in $[2b(0), 2b(0) + l]$ where the algorithm can put point mass.
		 
		  For the ``$\lambda$-problems'' \eqref{PROGRAM_P_1_short} and \eqref{PROGRAM_D_2} we discretise the interval $(0, t_0]$ by choosing $n_\lambda=100$ equidistant points where the algorithm can put mass. Table \ref{TabNumRes} gives the optimal values of all these four programs.
		  
	Moreover,  for the functions 
	\begin{align*}
		r_{\mu_{1,n}}\big(t, b(t)\big) \coloneqq \int_{[2b(0), \infty)} r_\theta\big(t, b(t)\big) \mu_{1,n}(d\theta)
	\end{align*}
	 and $r_{\mu_{2,n}}\big(t, b(t)\big)$ with $t$ between $0$ and $t_0=1$ we consider the minima and maxima of $r_{\mu_{1,n}}^{-1}$ and $r_{\mu_{2,n}}^{-1}$. According to Proposition \ref{PropNumBound}, this will give us an idea of the quality of the approximation of the true distribution function $F$ with the numerical c.d.f.\ $\widetilde{F}$. The results can be found in Table \ref{TabR1R2}. Also note that the computational times for all algorithms were very fast (less than 1 sec). 
	
		\begin{table}[hpbt]
		\setlength{\tabcolsep}{8pt} 
		\renewcommand{\arraystretch}{1.3} 
		\begin{center}
			\begin{tabular}{ |c|c|c|c|c| }
				\hline
				$b = b(t)$ & $d_{1,n}$ & $p_{1,n_\lambda}$ & $d_{2,n_\lambda}$ & $p_{2,n}$ 
				\\
				\hline
				$1+t$ & 0.1353353 & 0.1353353 & 0.1353353 & 0.1353353
				\\
				$\sqrt{1 + t}$ & 0.1274203 & 0.1274203 & 0.1274203 & 0.1274203 
				\\
				$\log(2+t)$ & 0.2364878 & 0.2364878 & 0.2364878 & 0.2364878  
				\\
				$1+t^2$ & 0.1353353 & 0.1353353 & 0.9801987 & 0.9980020 
				\\
				\hline
			\end{tabular}
			\caption[Numerical results for one-sided boundaries]{\label{TabNumRes}Numerical results for the optimal values for \eqref{PROGRAM_D_1_short} with discretisations as defined above.}
		\end{center}
	\end{table}
	\begin{table}[hpbt]
		\begin{center}
			\setlength{\tabcolsep}{8pt} 
			\renewcommand{\arraystretch}{1.3} 
			\begin{tabular}{ |c|c|c|c|c| }
				\hline
				$b = b(t)$ & $\min r_{\mu_{1,n}}^{-1}$ & $\max r_{\mu_{1,n}}^{-1}$ & $\min r_{\mu_{2,n}}^{-1}$ & $\max r_{\mu_{2,n}}^{-1}$\\
				\hline
				$1+t$ & 0.999999999953 & 1.000000000016 & 0.999999999996 & 1.000000000037 \\
				$\sqrt{1 + t}$ & 0.999999986788 & 1.002635679587 & 0.992114401333 & 1.000003242113 \\
				$\log(2+t)$ & 0.999999954916 & 1.001276226376 & 0.998198226881 & 1.000000116299 \\
				$1+t^2$ & 1.000000005738 & 7.389040705848 & 0.135606226531 & 1.001999983886 \\
				\hline
			\end{tabular}
			\caption[Numerical results for $r_{\mu_{1,n}}^{-1}$ and $r_{\mu_{2,n}}^{-1}$]{\label{TabR1R2}Numerical results for the minima and maxima of $r_{\mu_{1,n}}^{-1}$ and $r_{\mu_{2,n}}^{-1}$ on the interval $(0, t_0]$.}
		\end{center}
	\end{table}
	
	We immediately see in Table \ref{TabNumRes} that in the case of the linear boundary $b(t) = 1+t$ the values of $d_{1,n}, p_{1,n_\lambda}, d_{2,n_\lambda}$ and $p_{2,n}$ agree in the first $8$ digits. So, we can heuristically confirm strong duality in both set-ups as well as that the conditions for representability from Theorem~\ref{ThmRepEqual} are met, i.e., $d_1 = p_2$ and $p_2 = d_1$. This is, of course, not surprising as we already know that linear boundaries $b$ are representable, recall Remark \ref{rem:linear}. Therefore, this case can serve as a check of our algorithm.
	
\begin{figure}[hpbt]
		\begin{center}
			\includegraphics[scale=0.5]{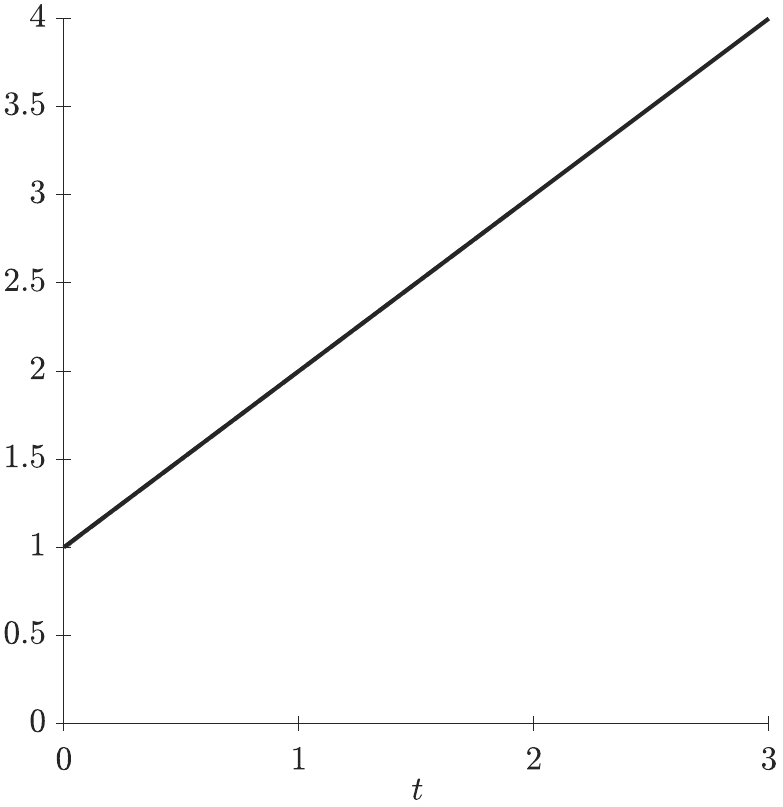}
			\includegraphics[scale=0.5]{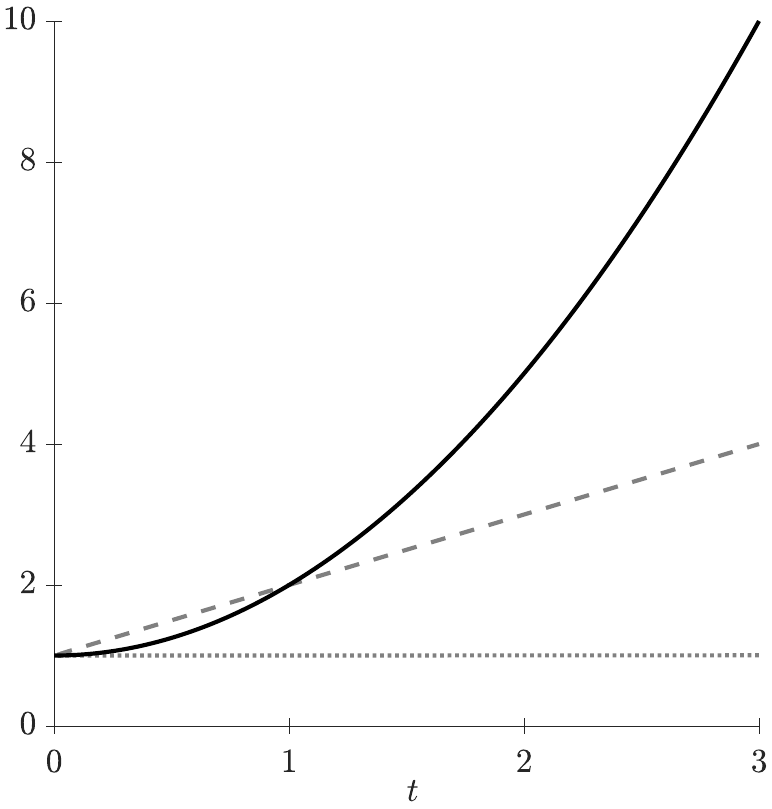}
			
			\includegraphics[scale=0.5]{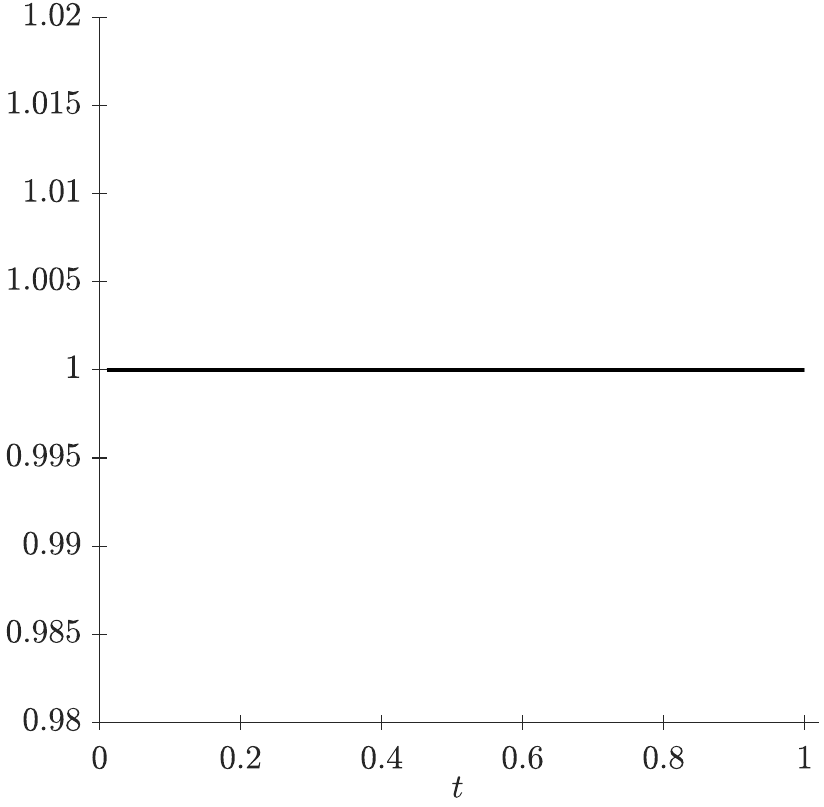}
			\includegraphics[scale=0.55]{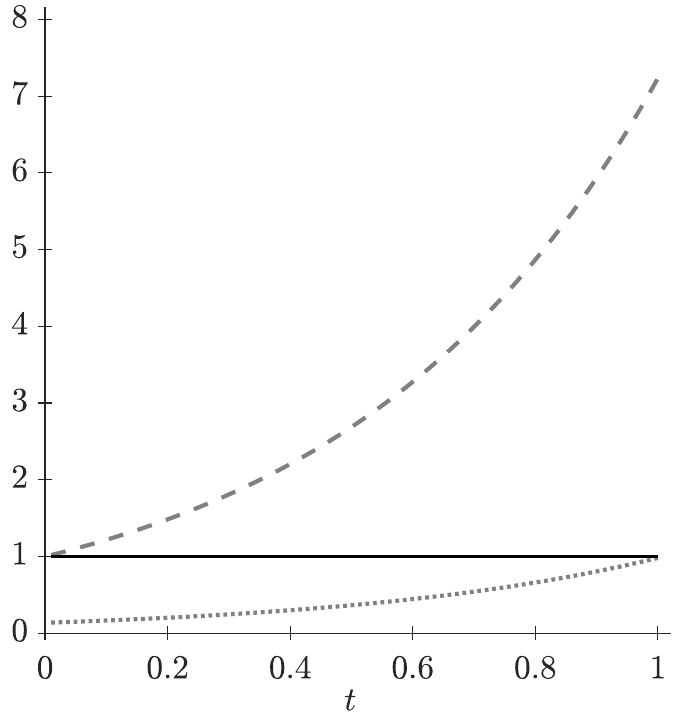}
			
			\includegraphics[scale=0.5]{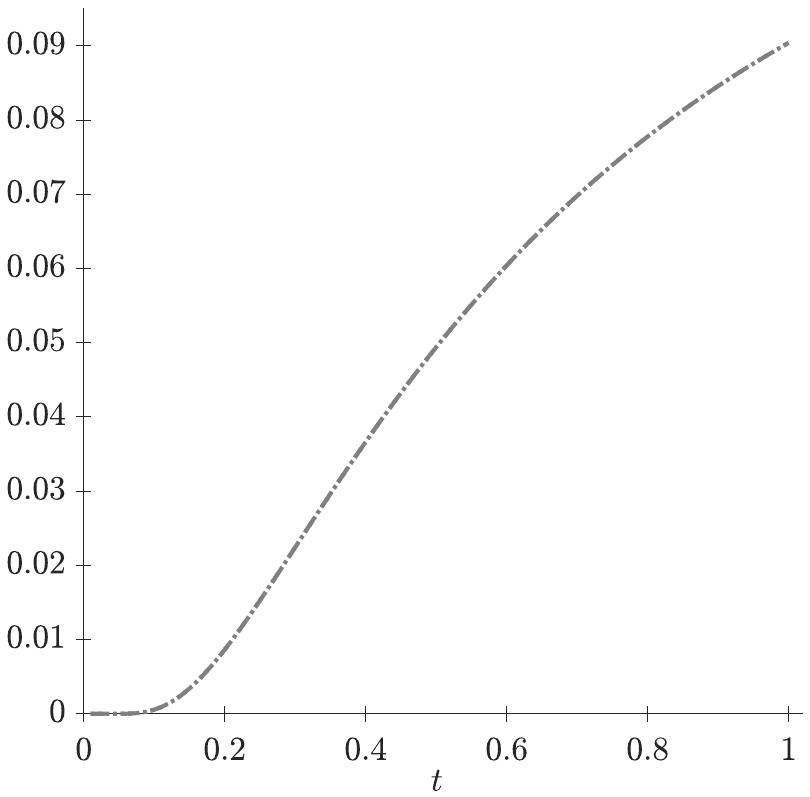}
			\includegraphics[scale=0.5]{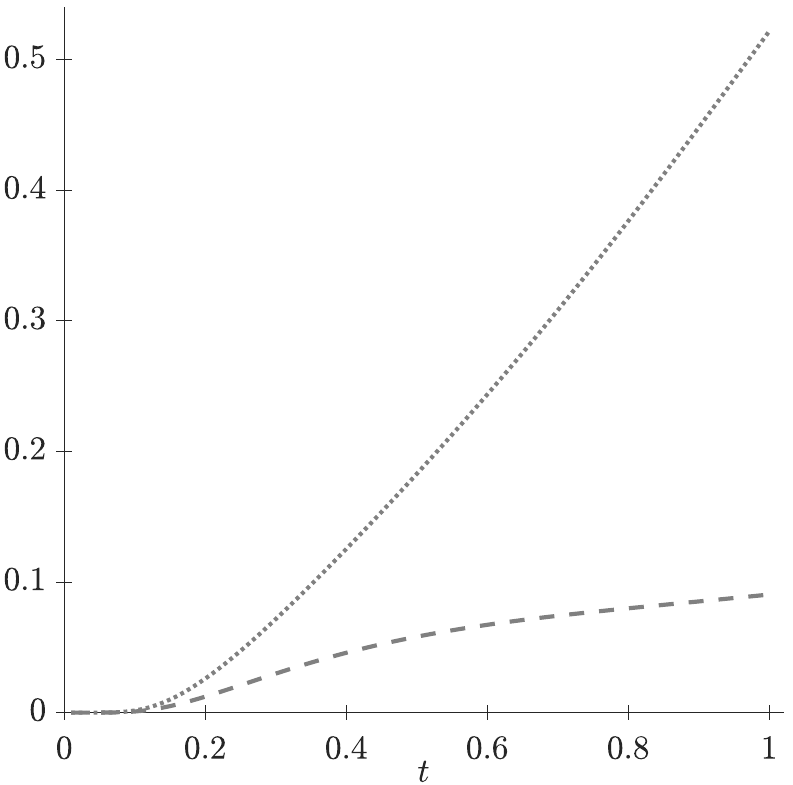}
			\caption[Numerical results for $b(t) = 1+t$ and $b(t)=1+t^2$]{Numerical results for $b(t) = 1+t$ on the left and for $b(t)=1+t^2$ on the right  for $t_0=1$.\\
				\emph{Top}: graphs of the original boundary (black) which overlaps the boundaries generated by $\mu_{1,n}$ (dashed) and $\mu_{2,n}$ (dotted).\\
				\emph{Middle}: corresponding values of $r_{\mu_{1,n}}\big(t, b(t)\big)$ (dashed) and $r_{\mu_{2,n}}\big(t, b(t)\big)$ (dotted) for $t\in[0,1]$ overlapped by the function constant $1$ (black).\\
				\emph{Bottom}: numerical approximations of the distribution functions $F_{\mu_1}$ (dashed) and $F_{\mu_2}$ (dotted).}
			\label{GraphsLinear}
		\end{center}
	\end{figure}

\begin{figure}[hpbt]
		\begin{center}
			\includegraphics[scale=0.5]{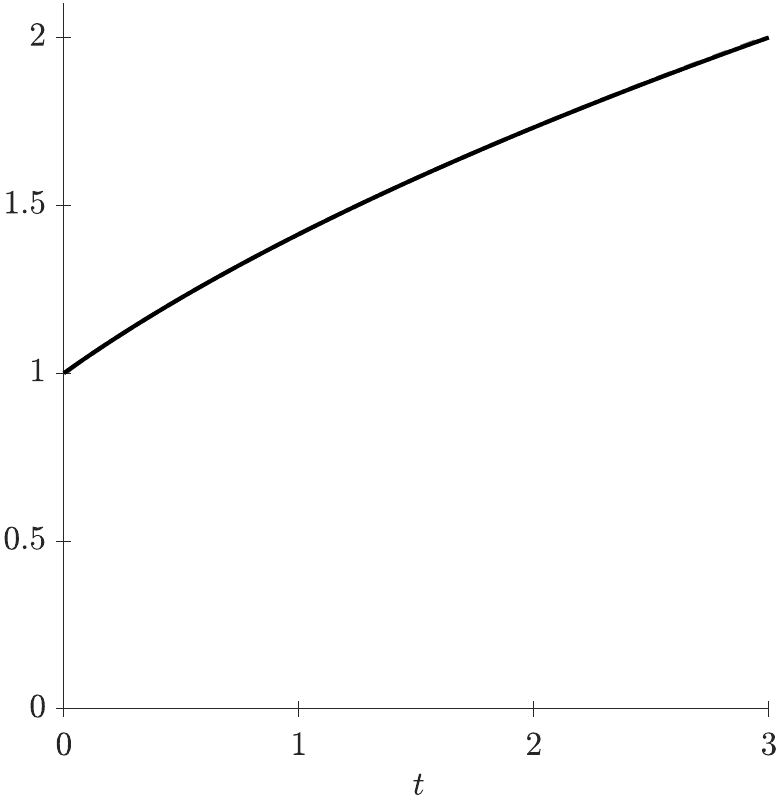}\hspace{0.98cm}
			\includegraphics[scale=0.5]{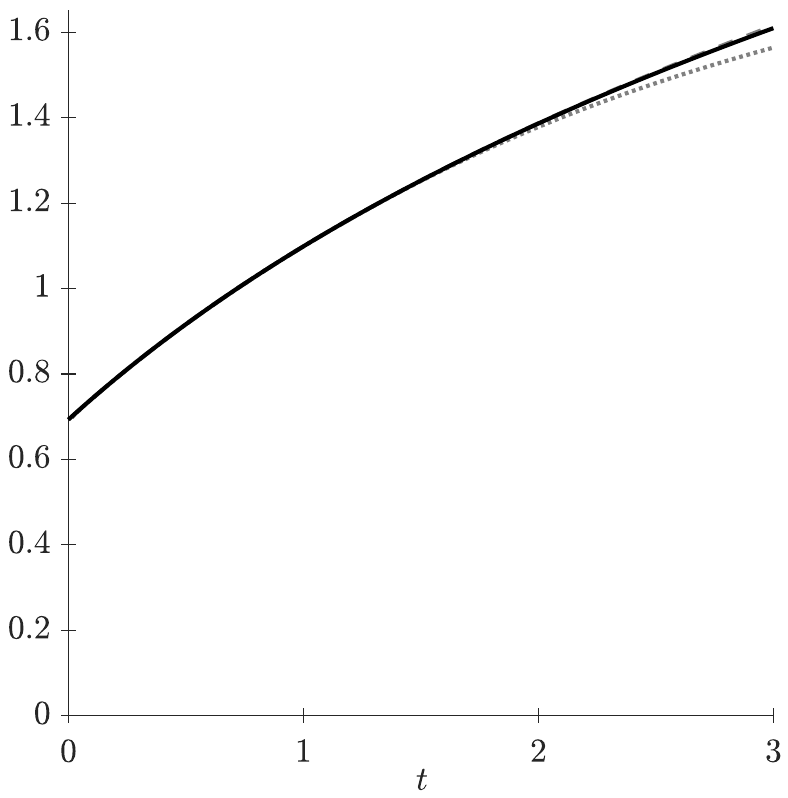}
			
			\includegraphics[scale=0.5]{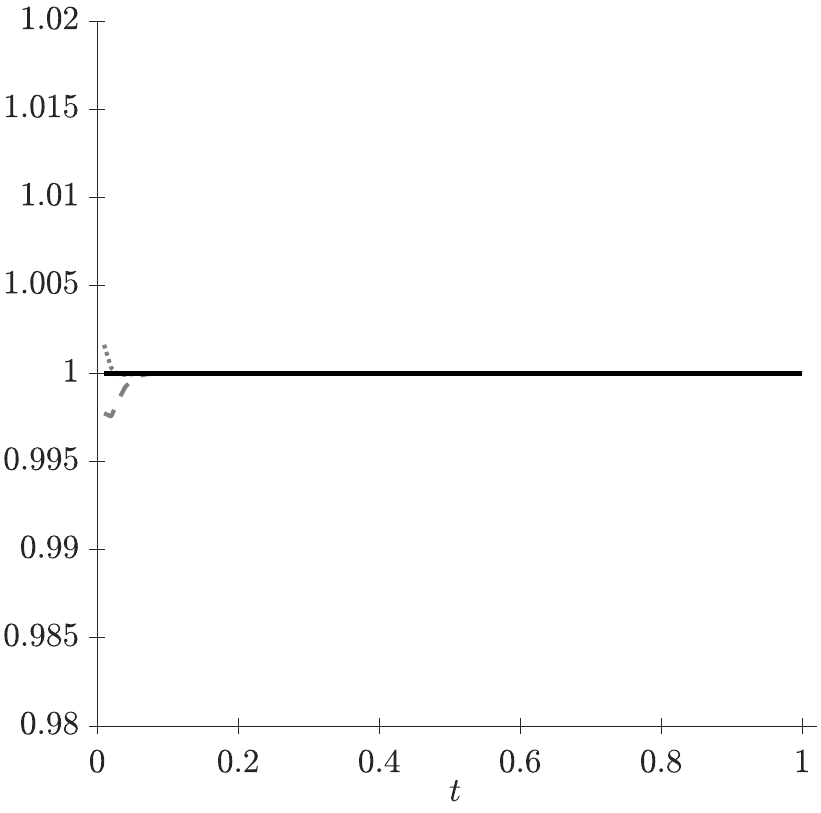}\hspace{0.5cm}
			\includegraphics[scale=0.51]{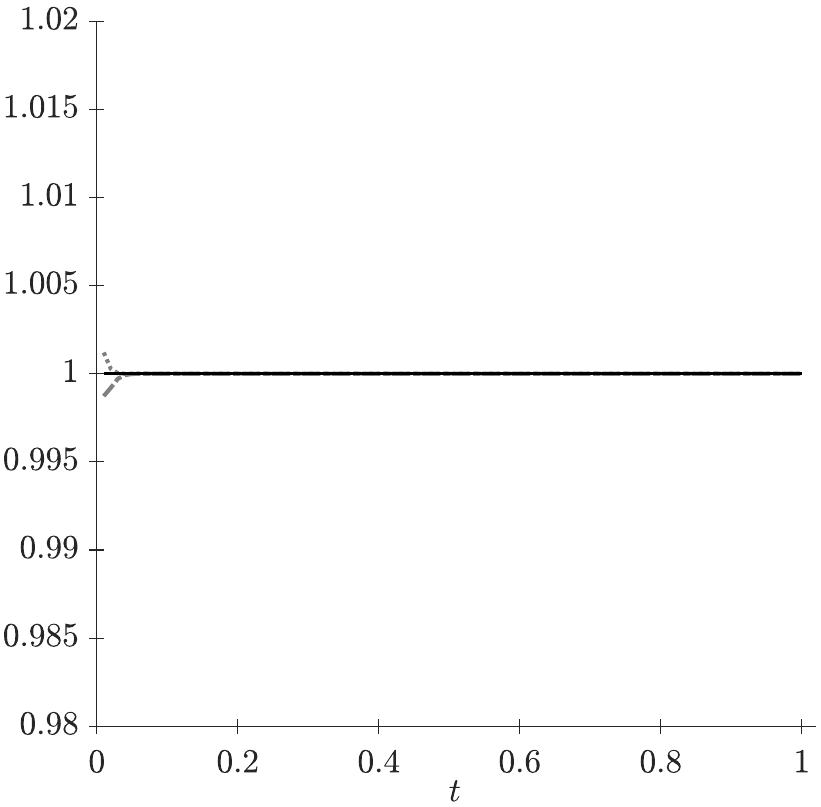}
		
			\includegraphics[scale=0.5]{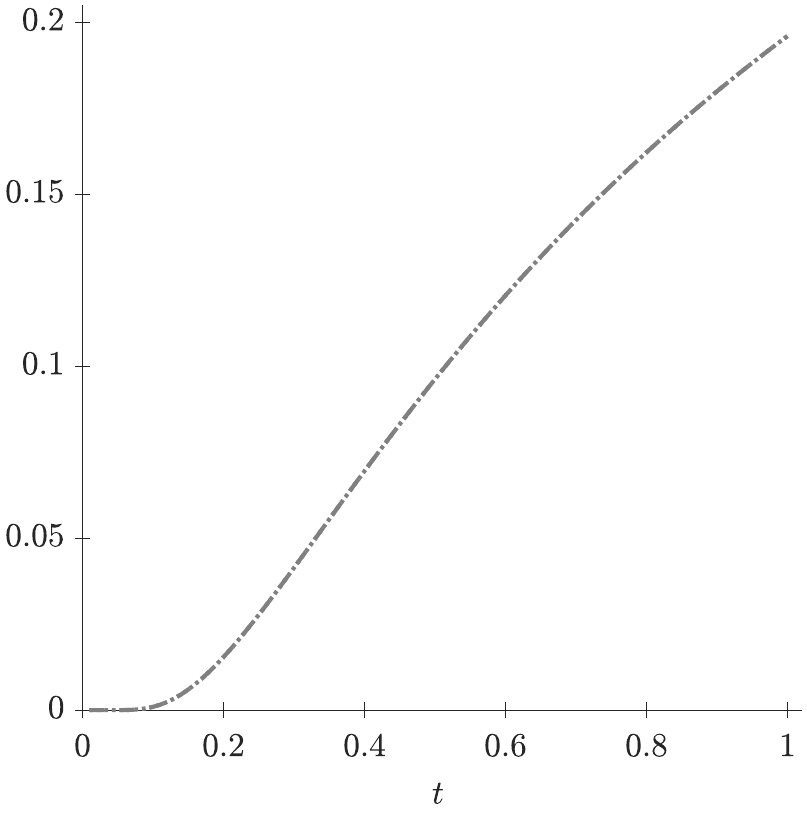}\hspace{0.73cm}
			\includegraphics[scale=0.5]{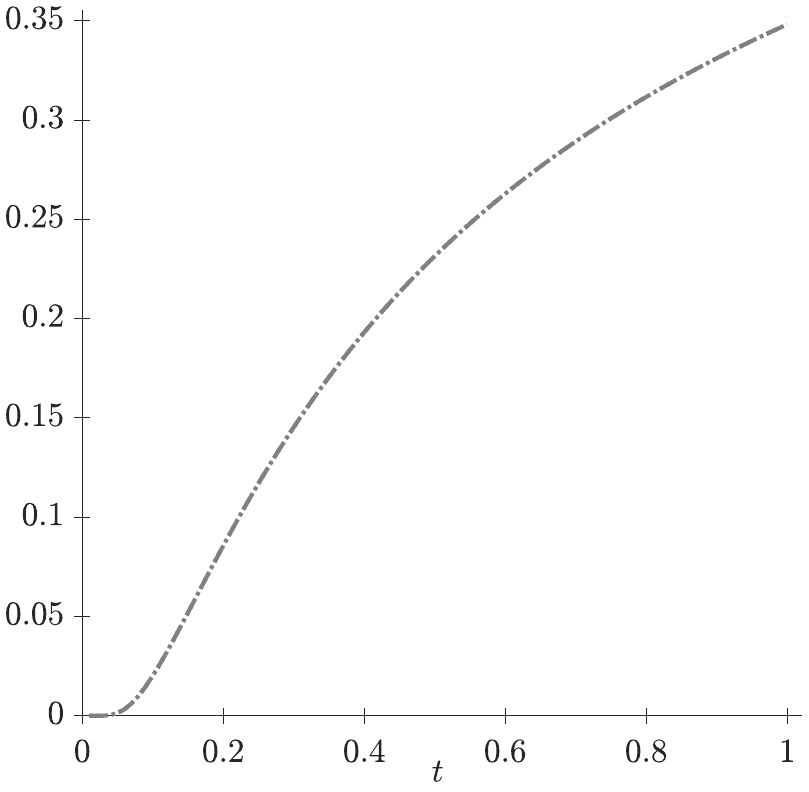}
			\caption[Numerical results for $b(t) = \sqrt{1+t}$ and $b(t) = \log(2+t)$]{Numerical results for $b(t) = \sqrt{1+t}$ on the left and for $b(t) = \log(2+t)$ on the right for $t_0=1$.\\
				\emph{Top}: graphs of the original boundary (black) which overlaps the boundaries generated by $\mu_{1,n}$ (dashed) and $\mu_{2,n}$ (dotted). The boundaries coincide even beyond the controlled interval $(0,1]$.\\
				\emph{Middle}: corresponding values of $r_{\mu_{1,n}}\big(t, b(t)\big)$ (dashed) and $r_{\mu_{2,n}}\big(t, b(t)\big)$ (dotted) for $t\in[0,1]$ compared with the function constant $1$ (black).\\
				\emph{Bottom}: numerical approximations of the distribution functions $F_{\mu_1}$ (dashed) and $F_{\mu_2}$ (dotted).}
			\label{GraphsSqrt}
		\end{center}
	\end{figure}

	Moreover, we consider the boundaries $b_{\mu_{1,n}}$ and $b_{\mu_{2,n}}$ generated by the numerical solutions $\mu_{1,n}$ and $\mu_{2,n}$ and compare these boundaries to the boundary $b$ that was the input for the optimisation problems. For our four input boundaries Figure \ref{GraphsLinear} and \ref{GraphsSqrt} depict the graphs of these three boundaries in the top figure, show the value of $r_{\mu_{1,n}}\big(t, b(t)\big)$ and $r_{\mu_{2,n}}\big(t, b(t)\big)$ between $0$ and $1$ in the middle figure and in the bottom the distribution function obtained by substituting the representing measure $\mu$ with the numerical solutions $\mu_{1,n}$ and $\mu_{2,n}$, respectively, i.e.,
	\begin{align*}
		F_{\mu_i}(t) \coloneqq 1-\Phi \left( \frac{b(t)}{\sqrt{t}} \right) + \int_{(0,\infty)}  \Phi \left( \frac{b(t)-\theta}{\sqrt{t}} \right) \mu_{i,n}(d \theta)
	\end{align*}
	for $i = 1,2$. Note that $F_{\mu_i}$ yields an approximation of the true distribution function $F$ in analytical form.
	We can see on the left of Figure \ref{GraphsLinear} that for the linear boundary $b(t)=1+t$  both numerical boundaries perfectly replicate the original boundary. In particular, the graphs of $r_{\mu_{i,n}}\big(t, b(t)\big), i=1, 2,$ which are equal to $1$ support that $b$ is represented by $\mu_{i,n}$. This is not very surprising as it is well-known that linear boundaries are representable by measures $\mu$ which only put mass into $2b(0)$ which both $\mu_{1,n}$ and $\mu_{2,n}$ do. Moreover, we can see in Table \ref{TabR1R2} that $r_{\mu_{1,n}}^{-1}$ and $r_{\mu_{2,n}}^{-1}$ deviate from $1$ by less than $10^{-10}$. The deviation should be $0$ and can probably be attributed to small numerical rounding errors. In particular, we know by Proposition \ref{PropNumBound} that we get a very good approximation of the distribution function $F$.
	
	Let us now consider the boundaries $b(t) = \sqrt{1+t}$ and $b(t) = \log(2+t)$, i.e., boundaries which are concave and monotone increasing. For both boundaries, we  observe (compare Table \ref{TabNumRes}) that the values of $d_{1,n}, p_{1,n_\lambda}, d_{2,n_\lambda}$ and $p_{2,n}$ agree in the first $8$ digits, i.e., we  again numerically confirm both strong duality in both set-ups as well as the conditions for representability from Theorem~\ref{ThmRepEqual}.
	Even though we choose $t_0=1$ we observe in the top of Figure \ref{GraphsSqrt} that the boundaries are very well replicated by the boundaries generated by $\mu_{1,n}$ and $\mu_{2,n}$ up to $t=2$ with just slight deviations between $t=2$ and $t=3$. The graphs for $r_{\mu_{1,n}}\big(t, b(t)\big)$ and $r_{\mu_{2,n}}\big(t, b(t)\big)$ are virtually indistinguishable from $1$ and we can see in Table \ref{TabR1R2} that $r_{\mu_{1,n}}^{-1}$ and $r_{\mu_{2,n}}^{-1}$ deviate from $1$ by less than $10^{-2}$. Thus, the graphs of the numerical distribution functions $F_{\mu_i}, i=1,2$, are very exact approximations for $F$ due to Proposition \ref{PropNumBound}. Moreover, note that Proposition~\ref{PropNumBound} gives a rather conservative estimate of the approximation error, so the true approximation error is probably much less.

	Now we turn our attention to $b(t) = 1+t^2$. We immediately see from Table \ref{TabNumRes} that $d_{1,n}$ and $p_{1,n_\lambda}$ agree in the first $8$ digits and $d_{2,n_\lambda}$ and $p_{2,n}$ agree in the first $2$ digits, so we can safely assume that strong duality holds in both set-ups. However, the gaps between the two different set-ups are very large.

	 We immediately see on the right of Figure \ref{GraphsLinear} why the algorithms does not produce sensible results: both programs only allow for concave boundaries and the ``most convex'' the program can do is a linear boundary. Unsurprisingly, $r_{\mu_{1,n}}\big(t, b(t)\big)$ and $r_{\mu_{2,n}}\big(t, b(t)\big)$ are for the most part very far away from $1$. 
  Expectedly, the numerical versions of the distribution functions are vastly different from one another.
	
	Finally, we investigate whether the conditions of Theorem \ref{ThmRepSeq2} are numerically met. Whether $\lambda_{1,n_\lambda}$ does not put mass into $t_0$ is hard to verify numerically, since the mass could be so small that it numerically vanishes. But we can investigate the mass which $\lambda_{2,n_\lambda}$ puts close to $t_0$ for increasing $n_\lambda$. Due to the discretisation of the $t$ axis, we look for mass in the interval $((n-1)t_0/n, t_0]$. In order to make sure that the mass near $t_0$ does not vanish, we need to ensure that $x_0$ is not too far below $b(t_0)$ as the {hitting probability} might become $0$ numerically. So, we choose $x_0 = b(t_0) -0.1$. All other parameters are chosen in the same way as before. Linear $b$ are omitted in this test, as we already know that they are representable and  $\bar{\lambda}$, the conditional last passage time distribution of the Brownian motion $W$ to $b$, can be shown to be an optimal measure in this case, cf.\ \cite[Corollary 2.16]{DissOskar}.  Table \ref{TableLambda} shows the values we obtained. 
	\begin{table}[hpbt]
		\begin{center}
			\setlength{\tabcolsep}{8pt} 
			\renewcommand{\arraystretch}{1.2} 
			\begin{tabular}{|c|c|c|c|}
				\hline
				$b = b(t)$ & $n_\lambda = 100$ & $n_\lambda = 200$ & $n_\lambda = 500$ 
				\\
				\hline
				$\sqrt{1 + t}$ & $0.233$ & $0.282$ & $0.276$  
				\\
				$\log(2+t)$ & $0.293$ & $0.290$ & $0.296$ 
				\\
				$1+t^2$ & $0.000$ & $0.000$ & $0.000$ 
				\\
				\hline
			\end{tabular}
			\caption[Numerical results for $\lambda_{2,n_\lambda}((n_\lambda-1)t_0/n_\lambda, t_0\rbrack$]{\label{TableLambda}Numerical results for $\lambda_{2,n_\lambda}\big((n_\lambda-1)t_0/n_\lambda, t_0\big]$. Values are rounded to 3 decimal points.}
		\end{center}
	\end{table}
	
	We immediately see in Table \ref{TableLambda} that for the concave boundaries the measure $\lambda_{2,n}$ always puts mass near $t_0$ while we increase $n_\lambda$. This is yet another very strong indicator that concave boundaries are representable according to Theorem \ref{ThmRepSeq2}. In contrast, for the convex boundary there is no mass near $t_0$ for any choice of $n_\lambda$, which is an indicator that the condition from Theorem \ref{ThmRepSeq2} is not only sufficient but might also be necessary.
	
	To sum up, the algorithm shows perfect replication of linear boundaries and very good replication of concave boundaries even beyond the controlled interval $(0, t_0]$.

	\section{Extension to two-sided boundaries}\label{ChapTwoSidedBoundaries}

	Finally, we shortly comment on the (inverse) method of images for two-sided boundaries. For more details we refer to \cite{DissOskar}. 
	
	Denote by $b_1$ the lower bound and by $b_2$ the upper bound for which we are interested in the first passage time distribution of a one-dimensional Brownian motion $W$. In particular, assume 
	\begin{itemize}
		\item $b_1\leq b_2$,
		\item $b_1(0) < W_0 < b_2(0)$
	\end{itemize}
	and define 
	\begin{align*}
		\tau =\inf\left\{t\in[0,\infty)\colon W_t \notin \big(b_1(t), b_2(t)\big)\right\}.
	\end{align*}

 An analogue to Proposition \ref{prop:meth_images} for the problem with two-sided boundaries can be proven which allows to formulate two linear programs which approximate the measure $\mu$ representing $b_1$ and $b_2$ from below and from above, respectively, and these programs are analogous to the problems for the classic method of images with one boundary.   
In particular, the linear programs and their formal duals can implement in the same way as in Section \ref{SecNumStuRep}, cf.\ \cite[Section 3.4]{DissOskar}.


\small

	
	\begin{thebibliography}{10}
		
		\bibitem{AB06}
		C.~D. Aliprantis and K.~C. Border.
		\newblock {\em Infinite {D}imensional {A}nalysis}.
		\newblock Springer, Berlin, third edition, 2006.
		
		\bibitem{BAC00}
		L.~Bachelier.
		\newblock {\em Th{\'e}orie de la sp{\'e}culation}.
		\newblock Gauthier-Villars, 1900.
		
		\bibitem{BAU01}
		H.~Bauer.
		\newblock {\em Measure and {I}ntegration {T}heory}, volume~26 of {\em De
			Gruyter Studies in Mathematics}.
		\newblock Walter de Gruyter \& Co., Berlin, 2001.
		
		\bibitem{BN05}
		K.~Borovkov and A.~Novikov.
		\newblock Explicit bounds for approximation rates of boundary crossing
		probabilities for the {W}iener process.
		\newblock {\em J. Appl. Probab.}, 42(1):82--92, 2005.
		
		\bibitem{BRE67}
		L.~Breiman.
		\newblock First exit times from a square root boundary.
		\newblock In {\em Proc. {F}ifth {B}erkeley {S}ympos. {M}ath. {S}tatist. and
			{P}robability ({B}erkeley, {C}alif., 1965/66), {V}ol. {II}: {C}ontributions
			to {P}robability {T}heory, {P}art 2}, pages 9--16. Univ. California Press,
		Berkeley, Calif., 1967.
		
		\bibitem{CHE06}
		L.~Cheng, X.~Chen, J.~Chadam, and D.~Saunders.
		\newblock Analysis of an inverse first passage problem from risk management.
		\newblock {\em SIAM J. Math. Anal.}, 38(3):845--873, 2006.
		
		\bibitem{CHR14}
		S.~Christensen.
		\newblock A method for pricing {A}merican options using semi-infinite linear
		programming.
		\newblock {\em Math. Finance}, 24(1):156--172, 2014.
		
		\bibitem{MR4635691}
		S.~Christensen, S.~Fischer, and O.~Hallmann.
		\newblock Uniqueness of first passage time distributions via {F}redholm
		integral equations.
		\newblock {\em Statist. Probab. Lett.}, 203:Paper No. 109912, 8, 2023.
		
		\bibitem{CKL22}
		S.~Christensen, J.~Kallsen, and M.~Lenga.
		\newblock Are {A}merican options {E}uropean after all?
		\newblock {\em Ann. Appl. Probab.}, 32(2):853--892, 2022.
		
		\bibitem{CGJ08}
		D.~Coculescu, H.~Geman, and M.~Jeanblanc.
		\newblock Valuation of default-sensitive claims under imperfect information.
		\newblock {\em Finance Stoch.}, 12(2):195--218, 2008.
		
		\bibitem{DAN82}
		H.~E. Daniels.
		\newblock Sequential tests constructed from images.
		\newblock {\em Ann. Statist.}, 10(2):394 -- 400, 1982.
		
		\bibitem{DAN00}
		H.~E. Daniels.
		\newblock The first crossing-time density for {B}rownian motion with a
		perturbed linear boundary.
		\newblock {\em Bernoulli}, 6(4):571--580, 2000.
		
		\bibitem{NAR01}
		E.~Di~Nardo, A.~G. Nobile, E.~Pirozzi, and L.~M. Ricciardi.
		\newblock A computational approach to first-passage-time problems for
		{G}auss-{M}arkov processes.
		\newblock {\em Adv. in Appl. Probab.}, 33(2):453--482, 2001.
		
		\bibitem{DUR71}
		J.~Durbin.
		\newblock Boundary-crossing probabilities for the {B}rownian motion and
		{P}oisson processes and techniques for computing the power of the
		{K}olmogorov-{S}mirnov test.
		\newblock {\em J. Appl. Probability}, 8:431--453, 1971.
		
		\bibitem{DUR85}
		J.~Durbin.
		\newblock The first-passage density of a continuous {G}aussian process to a
		general boundary.
		\newblock {\em J. Appl. Probab.}, 22(1):99--122, 1985.
		
		\bibitem{DUR92}
		J.~Durbin.
		\newblock The first-passage density of the {B}rownian motion process to a
		curved boundary.
		\newblock {\em J. Appl. Probab.}, 29(2):291--304, 1992.
		
		\bibitem{FER82B}
		B.~Ferebee.
		\newblock Tests with parabolic boundary for the drift of a {W}iener process.
		\newblock {\em Ann. Statist.}, 10(3):882 -- 894, 1982.
		
		\bibitem{FB09}
		E.~Freitag and R.~Busam.
		\newblock {\em Complex analysis}.
		\newblock Universitext. Springer-Verlag, Berlin, second edition, 2009.
		
		\bibitem{GY96}
		H.~Geman and M.~Yor.
		\newblock Pricing and hedging double-barrier options: A probabilistic approach.
		\newblock {\em Math. Finance}, 6(4):365--378, 1996.
		
		\bibitem{GRO89}
		P.~Groeneboom.
		\newblock Brownian motion with a parabolic drift and {A}iry functions.
		\newblock {\em Probab. Theory Related Fields}, 81(1):79--109, 1989.
		
		\bibitem{DissOskar}
		O.~F. Hallmann.
		\newblock {\em Contributions to {F}irst {P}assage {T}ime {P}roblems of
			{B}rownian {M}otion}.
		\newblock PhD thesis, 2023.
		\newblock \url{https://macau.uni-kiel.de/receive/macau_mods_00003885}.
		
		\bibitem{HW01}
		J.~Hull and A.~White.
		\newblock Valuing credit default swaps {II}: Modeling default correlations.
		\newblock {\em J. Deriv.}, 8, 02 2001.
		
		\bibitem{IW10}
		S.~Ito, S.~Wu, T.~Shiu, and K.~Teo.
		\newblock A numerical approach to infinite-dimensional linear programming in
		{$L_1$} spaces.
		\newblock {\em J. Ind. Manag. Optim.}, 6:15--28, 2010.
		
		\bibitem{JKV09}
		S.~Jaimungal, A.~Kreinin, and A.~Valov.
		\newblock Integral equations and the first passage time of {B}rownian motions.
		\newblock \url{https://arxiv.org/abs/0902.2569}, preprint, 2009.
		
		\bibitem{JW17}
		Z.~Jin and L.~Wang.
		\newblock First passage time for {B}rownian motion and piecewise linear
		boundaries.
		\newblock {\em Methodol. Comput. Appl. Probab.}, 19(1):237--253, 2017.
		
		\bibitem{KAH08}
		N.~Kahale.
		\newblock Analytic crossing probabilities for certain barriers by {B}rownian
		motion.
		\newblock {\em Ann. Appl. Probab.}, 18(4):1424--1440, 2008.
		
		\bibitem{KHI33}
		A.~Khintchine.
		\newblock {\em Asymptotische Gesetze der Wahrscheinlichkeitsrechnung}, volume~1
		of {\em Ergebnisse der Mathematik und ihrer Grenzgebiete. 1. Folge}.
		\newblock Springer Berlin, Heidelberg, 1933.
		
		\bibitem{KLE14}
		A.~Klenke.
		\newblock {\em Probability theory}.
		\newblock Universitext. Springer, London, second edition, 2014.
		
		\bibitem{KI92}
		N.~Kunitomo and M.~Ikeda.
		\newblock Pricing options with curved boundaries.
		\newblock {\em Math. Finance}, 2(4):275--298, 1992.
		
		\bibitem{LW92}
		H.~C. Lai and S.~Y. Wu.
		\newblock Extremal points and optimal solutions for general capacity problems.
		\newblock {\em Math. Programming}, 54(1):87--113, 1992.
		
		\bibitem{LAI01}
		T.~L. Lai.
		\newblock Sequential analysis: some classical problems and new challenges.
		\newblock {\em Statist. Sinica}, 11(2):303--408, 2001.
		
		\bibitem{LAN93}
		S.~Lang.
		\newblock {\em Real and functional analysis}, volume 142 of {\em Graduate Texts
			in Mathematics}.
		\newblock Springer-Verlag, New York, third edition, 1993.
		
		\bibitem{LEN17}
		M.~Lenga.
		\newblock {\em Representable Options}.
		\newblock PhD thesis, 2017.
		\newblock \url{https://macau.uni-kiel.de/receive/diss_mods_00021000}.
		
		\bibitem{LER86}
		H.~R. Lerche.
		\newblock {\em Boundary crossing of {B}rownian motion}, volume~40 of {\em
			Lecture Notes in Statistics}.
		\newblock Springer-Verlag, Berlin, 1986.
		
		\bibitem{levy1965processus}
		P.~L{\'e}vy.
		\newblock Processus stochastiques et mouvement brownien, gauthier-villars,
		paris.
		\newblock {\em Reprinted by Editions J. Gabay, Paris}, 1965.
		
		\bibitem{LRD02}
		V.~S.~F. Lo, G.~O. Roberts, and H.~E. Daniels.
		\newblock Inverse method of images.
		\newblock {\em Bernoulli}, 8(1):53--80, 2002.
		
		\bibitem{MV92}
		R.~Meise and D.~Vogt.
		\newblock {\em Introduction to functional analysis}, volume~2 of {\em Oxford
			Graduate Texts in Mathematics}.
		\newblock The Clarendon Press, Oxford University Press, New York, 1997.
		
		\bibitem{NFK99}
		A.~Novikov, V.~Frishling, and N.~Kordzakhia.
		\newblock Approximations of boundary crossing probabilities for a {B}rownian
		motion.
		\newblock {\em J. Appl. Probab.}, 36(4):1019--1030, 1999.
		
		\bibitem{NOV81}
		A.~A. Novikov.
		\newblock A martingale approach to first passage problems and a new condition
		for {W}ald's identity.
		\newblock In M.~Arat{\'o}, D.~Vermes, and A.~V. Balakrishnan, editors, {\em
			{S}tochastic {D}ifferential {S}ystems}, pages 146--156, Berlin, Heidelberg,
		1981. Springer Berlin Heidelberg.
		
		\bibitem{PP74}
		C.~Park and S.~R. Paranjape.
		\newblock Probabilities of {W}iener paths crossing differentiable curves.
		\newblock {\em Pacific J. Math.}, 53:579--583, 1974.
		
		\bibitem{PES02}
		G.~Peskir.
		\newblock On integral equations arising in the first-passage problem for
		{B}rownian motion.
		\newblock {\em J. Integral Equations Appl.}, 14(4):397--423, 2002.
		
		\bibitem{Poe12}
		K.~P\"{o}tzelberger.
		\newblock Improving the {M}onte {C}arlo estimation of boundary crossing
		probabilities by control variables.
		\newblock {\em Monte Carlo Methods Appl.}, 18(4):353--377, 2012.
		
		\bibitem{PW01}
		K.~P\"{o}tzelberger and L.~Wang.
		\newblock Boundary crossing probability for {B}rownian motion.
		\newblock {\em J. Appl. Probab.}, 38(1):152--164, 2001.
		
		\bibitem{RED01}
		S.~Redner.
		\newblock {\em A {G}uide to {F}irst-{P}assage {P}rocesses}.
		\newblock Cambridge University Press, 2001.
		
		\bibitem{RED22}
		S.~Redner.
		\newblock A first look at first-passage processes.
		\newblock {\em Phys. A}, 631:Paper No. 128545, 2023.
		
		\bibitem{RY99}
		D.~Revuz and M.~Yor.
		\newblock {\em Continuous martingales and {B}rownian motion}, volume 293 of
		{\em Grundlehren der mathematischen Wissenschaften [Fundamental Principles of
			Mathematical Sciences]}.
		\newblock Springer-Verlag, Berlin, third edition, 1999.
		
		\bibitem{RSS84}
		L.~M. Ricciardi, L.~Sacerdote, and S.~Sato.
		\newblock On an integral equation for first-passage-time probability densities.
		\newblock {\em J. Appl. Probab.}, 21(2):302--314, 1984.
		
		\bibitem{RS73}
		H.~Robbins and D.~Siegmund.
		\newblock The expected sample size of some tests of power one.
		\newblock {\em Ann. Statist.}, 2(3):415 -- 436, 1974.
		
		\bibitem{robbins1985statistical}
		H.~Robbins and D.~Siegmund.
		\newblock Statistical tests of power one and the integral representation of
		solutions of certain partial differential equations.
		\newblock {\em Herbert Robbins Selected Papers}, pages 299--326, 1985.
		
		\bibitem{RS97}
		G.~O. Roberts and C.~F. Shortland.
		\newblock Pricing barrier options with time--dependent coefficients.
		\newblock {\em Math. Finance}, 7(1):83--93, 1997.
		
		\bibitem{ROC70}
		R.~T. Rockafellar.
		\newblock {\em Convex analysis}.
		\newblock Princeton Mathematical Series, No. 28. Princeton University Press,
		Princeton, N.J., 1970.
		
		\bibitem{ROC74}
		R.~T. Rockafellar.
		\newblock {\em Conjugate duality and optimization}.
		\newblock Society for Industrial and Applied Mathematics, Philadelphia, Pa.,
		1974.
		
		\bibitem{SAL88}
		P.~Salminen.
		\newblock On the first hitting time and the last exit time for a {B}rownian
		motion to/from a moving boundary.
		\newblock {\em Adv. in Appl. Probab.}, 20(2):411--426, 1988.
		
		\bibitem{SCH15}
		E.~Schr{\"o}dinger.
		\newblock Zur {T}heorie der {F}all- und {S}teigversuche an {T}eilchen mit
		{B}rownscher {B}ewegung.
		\newblock {\em Physik. Z.}, 16(16):289--295, 1915.
		
		\bibitem{SHE67}
		L.~A. Shepp.
		\newblock A first passage problem for the {W}iener process.
		\newblock {\em Ann. Math. Statist.}, 38:1912--1914, 1967.
		
		\bibitem{SIE86}
		D.~Siegmund.
		\newblock Boundary crossing probabilities and statistical applications.
		\newblock {\em Ann. Statist.}, 14(2):361 -- 404, 1986.
		
		\bibitem{SMI72}
		C.~S. Smith.
		\newblock A note on boundary-crossing probabilities for the {B}rownian motion.
		\newblock {\em J. Appl. Probab.}, 9:857--861, 1972.
		
		\bibitem{WP97}
		L.~Wang and K.~P\"{o}tzelberger.
		\newblock Boundary crossing probability for {B}rownian motion and general
		boundaries.
		\newblock {\em J. Appl. Probab.}, 34(1):54--65, 1997.
		
		\bibitem{ZIP13}
		P.~Zipkin.
		\newblock Linear programming and the inverse method of images.
		\newblock {\em Ann. Oper. Res.}, 208:227--243, 2013.
		
	\end{thebibliography}
\end{document}